\documentclass[12pt]{amsart}
\usepackage{amsmath}
\usepackage{amssymb}

\newtheorem{propo}{Proposition}[section]

\newtheorem{theor}[propo]{Theorem}
\newtheorem{lemma}[propo]{Lemma}
\theoremstyle{definition}
\newtheorem{defin}[propo]{Definition}
\newtheorem{examp}[propo]{Example}

\theoremstyle{remark}
\newtheorem{remar}[propo]{Remark}

\numberwithin{equation}{section}

\newcommand{\al }{\alpha }
\newcommand{\Aut }{\mathrm{Aut}}

\newcommand{\cC }{\mathcal{C}}

\newcommand{\cm }{c}
\newcommand{\Cm }{C}
\newcommand{\Dih }{\mathbb{D}}
\newcommand{\End }{\mathrm{End}}

\newcommand{\Hom }{\mathrm{Hom}}
\newcommand{\id }{\mathrm{id}}
\newcommand{\ka }{\kappa }
\newcommand{\Lin }{\mathrm{span}}

\newcommand{\ndN }{\mathbb{N}}

\newcommand{\ndZ }{\mathbb{Z}}
\newcommand{\Ob }{\mathrm{Ob}}
\newcommand{\om }{\omega }
\newcommand{\SG}[1]{\mathbb{S}_{#1}}
\newcommand{\re }{^\mathrm{re}}
\newcommand{\rsC }{\mathcal{R}}

\newcommand{\Stab }{\mathrm{Hom}}
\newcommand{\s }{\sigma }
\newcommand{\tr }{\mathrm{tr}\,}
\newcommand{\rfl }{\rho }
\newcommand{\Wg }{\mathcal{W}}

\title{Weyl groupoids with at most three objects}

\author{M.~Cuntz}
\address{Michael Cuntz,
Fachbereich Mathematik,
Universit\"at Kaiserslau\-tern,
Postfach 3049,
D-67653 Kaiserslautern, Germany}
\email{cuntz@mathematik.uni-kl.de}

\author{I.~Heckenberger}
\thanks{I.H. is
supported by the German Research Foundation (DFG) via a Heisenberg
fellowship}
\address{Istv\'an Heckenberger, Mathematisches Institut,
Ludwig-Maximili\-ans-Universit\"at M\"unchen,
Theresienstr. 39,
D-80333 M\"unchen, Germany}
\email{i.heckenberger@googlemail.com}

\begin{document}

\begin{abstract}
We adapt the generalization of root systems of the second author and H.~Yamane
to the terminology of category theory.
We introduce Cartan schemes, associated root systems and Weyl groupoids.
After some preliminary general results, we completely classify all
finite Weyl groupoids with at most three objects.
The classification yields that there exist infinitely many ``standard'', but
only 9 ``exceptional'' cases.

Key words: Nichols algebra, reflection, root system, Weyl groupoid
\end{abstract}

\maketitle

\section{Introduction}

In the last decades, root systems and their generalizations
have continuously led to many remarkable new results. In most cases
the motivation was to understand the structure of some generalization
of Lie algebras, for example Kac-Moody algebras or Lie superalgebras.

Following the plan of Andruskiewitsch and Schneider for a
classification of pointed Hopf algebras
\cite{a-AndrSchn98},\cite{a-AndrSchn00}, a new type of root systems
emerged \cite{a-Heck06a}. 
These are fundamental invariants of Nichols algebras of diagonal type, and are
crucial for the full classification of finite-dimensional Nichols algebras of
diagonal type \cite{p-Heck06b}.
In \cite{a-HeckYam08} an axiomatic definition of a generalization of root
systems was introduced, based on the main properties of the root systems of
Nichols algebras of diagonal type.
The class of these root systems includes properly
the reduced root systems in the sense of Bourbaki \cite{b-BourLie4-6} and the
root systems of Kac-Moody algebras \cite{b-Kac90}, but contains many
exceptional cases. The reduced root systems of simple Lie superalgebras also
fit naturally into the new framework \cite{a-HeckYam08}.
The results in \cite{p-AHS08} indicate that a large class of
Nichols algebras of non-diagonal type, such as those of finite group
type, presumably admits a root system satisfying the axioms given in
\cite{a-HeckYam08}.

The main aspect of the novel root systems is,
that one starts with a family of Cartan
matrices instead of a single Cartan matrix. Consequently,
the symmetry object is not a group but a
groupoid, the so called \textit{Weyl groupoid}.

There have also been many efforts to find ``root systems'' for
complex reflection groups \cite{a-BMM99},
the main achievement being the cyclotomic Hecke algebras.
In a connected groupoid, if one fixes an object $a$, then the morphisms
from $a$ to $a$ form a group $\Hom(a)$ which is isomorphic for
all choices of $a$.
It is easy to see that any finite group appears as $\Hom(a)$ for some
finite Weyl groupoid. So in particular, it will be interesting to
investigate the root systems for which $\Hom(a)$ is a complex
reflection group.

Matsumoto's theorem holds for Coxeter groupoids \cite{a-HeckYam08} and
hence there will be many nice properties of Coxeter groups which may
be translated to this new situation. However, there are some important
differences, for example the exchange condition only holds in a weak
version.

In this article, we focus on various different aspects of Weyl groupoids.
To stress the naturality of the construction,
we introduce new concepts for the definitions using the language of category
theory. A Weyl groupoid
is based on a set of Cartan matrices $\cC$ called a {\it Cartan scheme}.
For such a scheme $\cC$, we define root systems of type $\cC$
and their Weyl groupoid.
We stay in the general setting and deduce many useful results
about root systems, Cartan schemes and Weyl groupoids, extending the analysis
in \cite{a-HeckYam08}. We discuss standard Cartan schemes and their Weyl
groupoids: Regardless of the number of objects, these are defined
with the help of a single Cartan matrix, and are closely related to the
crystallographic Weyl groups. Then decompositions of Cartan schemes, root
systems, and their
Weyl groupoids are investigated and characterized. In Prop.~\ref{pr:decomp}
we prove
that a finite root system
for a given Cartan scheme is reducible if and only if
the family of Cartan matrices
is simultaneously decomposable, or equivalently, if one of the Cartan
matrices of the family is decomposable.

Then we turn our attention to the case of finite root systems.
The main theorems merge to the result:
\begin{theor}
Let $\Wg$ be the Weyl groupoid of a finite irreducible connected root system
with at most $3$ objects.
Then one of the following holds:
\begin{enumerate}
\item $\Wg$ is standard, i.e. all Cartan matrices are equal.
\item $\Wg$ is one of $9$ exceptional Weyl groupoids.
\end{enumerate}
\end{theor}
The main tool in our proofs is Thm.~\ref{th:Coxgr}, the proof of which is
given in \cite{a-HeckYam08}. It states that $\Wg $ is
generated by reflections and Coxeter relations.
For details on which Cartan matrices actually yield standard Weyl
groupoids and a description of the exceptional cases, see the
theorems \ref{th:o=2}, \ref{th:r2o3}, \ref{th:o=3,r=3} and \ref{th:r>3o3}.
As a consequence of our classification, we conclude in
Rems.~\ref{re:WgNichols2} and \ref{re:WgNichols3} that there exist
root systems associated to some non-standard Cartan schemes
which cannot be obtained as a root system of a
finite-dimensional Nichols algebra of diagonal type.

There are many open questions left. It is conceivable that there are
only finitely many non-standard irreducible connected Weyl
groupoids for a fixed number of objects.
Notice that there are infinitely many standard irreducible connected
Weyl groupoids with two objects, but that all irreducible connected
Weyl groupoids with three objects have rank less or equal four.

We use the symbols $\ndN $ and $\ndN _0$ for the set of positive and
nonnegative integers, respectively.

We want to thank G.~Malle for providing us with
Ex.~\ref{ex:rcwithgivenHom},
and H.-J.~Schneider for many interesting
discussions on the subject and for his help in looking for a good terminology. 

\section{Cartan schemes, root systems, and their Weyl groupoids}

First the generalization of root systems given in \cite[Def.\,2]{a-HeckYam08}
is reformulated in terms of category theory.

Let $I$ be a non-empty finite set and
$\{\al _i\,|\,i\in I\}$ the standard basis of $\ndZ ^I$.
Recall from \cite[\S 1.1]{b-Kac90}
that a generalized Cartan matrix $\Cm =(\cm _{ij})_{i,j\in I}$
is a matrix in $\ndZ ^{I\times I}$ such that
\begin{enumerate}
  \item[(M1)] $\cm _{ii}=2$ and $\cm _{jk}\le 0$ for all $i,j,k\in I$ with
    $j\not=k$,
  \item[(M2)] if $i,j\in I$ and $\cm _{ij}=0$, then $\cm _{ji}=0$.
\end{enumerate}

\begin{defin} \label{de:CS}
  Let $A$ be a non-empty set, $\rfl _i : A \to A$ a map for all $i\in I$,
  and $\Cm ^a=(\cm ^a_{jk})_{j,k \in I}$ a generalized Cartan matrix
  in $\ndZ ^{I \times I}$ for all $a\in A$. The quadruple
  \[ \cC = \cC (I,A,(\rfl _i)_{i \in I}, (\Cm ^a)_{a \in A})\]
  is called a \textit{Cartan scheme} if
  \begin{enumerate}
  \item[(C1)] $\rfl _i^2 = \id$ for all $i \in I$,
  \item[(C2)] $\cm ^a_{ij} = \cm ^{\rfl _i(a)}_{ij}$ for all $a\in A$ and
    $i,j\in I$.
  \end{enumerate}
  We say that $\cC $ is \textit{connected}, if
  the group $\langle \rfl _i\,|\,i\in I\rangle \subset \Aut (A)$ acts
  transitively on $A$, that is, if
  for all $a,b\in A$ with $a\not=b$ there exist $n\in \ndN $,
  $a_1,a_2,\ldots ,a_n\in A$, and $i_1,i_2,\ldots ,i_{n-1}\in I$ such that
  \[a_1=a,\quad a_n=b,\quad a_{j+1}=\rfl _{i_j}(a_j)\quad
  \text{for all $j=1,\ldots ,n-1$.}\]
  Two Cartan schemes $\cC =\cC (I,A,(\rfl _i)_{i\in I},(\Cm ^a)_{a\in A})$
  and $\cC '=\cC '(I',A'$,
  $(\rfl '_i)_{i\in I'},({\Cm '}^a)_{a\in A'})$
  are termed
  \textit{equivalent}, if there are bijections $\varphi _0:I\to I'$
  and $\varphi _1:A\to A'$ such that
  \begin{align*}
    \varphi _1(\rfl _i(a))=\rfl '_{\varphi _0(i)}(\varphi _1(a)),
    \qquad
    \cm ^{\varphi _1(a)}_{\varphi _0(i) \varphi _0(j)}=\cm ^a_{i j}
  \end{align*}
  for all $i,j\in I$ and $a\in A$.

  Let $\cC = \cC (I,A,(\rfl _i)_{i \in I}, (\Cm ^a)_{a \in A})$ be a
  Cartan scheme. For all $i \in I$ and $a \in A$ define $\s _i^a \in
  \Aut(\ndZ ^I)$ by
  \begin{align}
    \s _i^a (\al _j) = \al _j - \cm _{ij}^a \al _i \qquad
    \text{for all $j \in I$.}
    \label{eq:sia}
  \end{align}
  The \textit{Weyl groupoid of} $\cC $
  is the category $\Wg (\cC )$ such that $\Ob (\Wg (\cC ))=A$ and
  the morphisms are generated by the maps
  $\s _i^a\in \Hom (a,\rfl _i(a))$ with $i\in I$, $a\in A$.
  Formally, for $a,b\in A$ the set $\Hom (a,b)$ consists of the triples
  $(b,f,a)$, where
  \[ f=\s _{i_n}^{\rfl _{i_{n-1}}\cdots \rfl _{i_1}(a)}\cdots
    \s _{i_2}^{\rfl _{i_1}(a)}\s _{i_1}^a \]
  and $b=\rfl _{i_n}\cdots \rfl _{i_2}\rfl _{i_1}(a)$ for some
  $n\in \ndN _0$ and $i_1,\ldots ,i_n\in I$.
  The composition is induced by the group structure of $\Aut (\ndZ ^I)$:
  \[ (a_3,f_2,a_2)\circ (a_2,f_1,a_1) = (a_3,f_2f_1, a_1)\]
  for all $(a_3,f_2,a_2),(a_2,f_1,a_1)\in \Hom (\Wg (\cC ))$.
  By abuse of notation we will write
  $f\in \Hom (a,b)$ instead of $(b,f,a)\in \Hom (a,b)$.
 
  The cardinality of $I$ is termed the \textit{rank of} $\Wg (\cC )$.
\end{defin}

The Weyl groupoid $\Wg (\cC )$ of a Cartan scheme $\cC $ is a groupoid.
Indeed, (M1) implies that
$\s _i^a\in \Aut (\ndZ ^I)$ is a reflection for all $i\in I$ and $a\in A$,
and hence the inverse of $\s _i^a\in \Hom (a,\rfl _i(a))$ is
$\s _i^{\rfl _i(a)}\in \Hom (\rfl _i(a),a)$
by (C1) and (C2). Therefore each morphism of $\Wg (\cC )$ is an isomorphism.

If $\cC $ and $\cC '$ are equivalent Cartan schemes, then $\Wg (\cC )$ and
$\Wg (\cC ')$ are isomorphic groupoids.

Recall that a groupoid $G$ is called \textit{connected},
if for each $a,b\in \Ob (G)$ the class $\Hom (a,b)$ is non-empty.
Hence $\Wg (\cC )$ is a connected groupoid if and only if $\cC $ is a
connected Cartan scheme.

\begin{defin} \label{de:RSC}
  Let $\cC =\cC (I,A,(\rfl _i)_{i\in I},(\Cm ^a)_{a\in A})$ be a Cartan
  scheme. For all $a\in A$ let $R^a\subset \ndZ ^I$, and define
  $m_{i,j}^a= |R^a \cap (\ndN _0 \al _i + \ndN _0 \al _j)|$ for all $i,j\in
  I$ and $a\in A$. We say that
  \[ \rsC = \rsC (\cC , (R^a)_{a\in A}) \]
  is a \textit{root system of type} $\cC $, if it satisfies the following
  axioms.
  \begin{enumerate}
    \item[(R1)]
      $R^a=R^a_+\cup - R^a_+$, where $R^a_+=R^a\cap \ndN _0^I$, for all
      $a\in A$.
    \item[(R2)]
      $R^a\cap \ndZ \al _i=\{\al _i,-\al _i\}$ for all $i\in I$, $a\in A$.
    \item[(R3)]
      $\s _i^a(R^a) = R^{\rfl _i(a)}$ for all $i\in I$, $a\in A$.
    \item[(R4)]
      If $i,j\in I$ and $a\in A$ such that $i\not=j$ and $m_{i,j}^a$ is
      finite, then
      $(\rfl _i\rfl _j)^{m_{i,j}^a}(a)=a$.
  \end{enumerate}
  If $\rsC $ is a root system of type $\cC $, then we say that
  $\Wg (\rsC )=\Wg (\cC )$ is the \textit{Weyl groupoid of} $\rsC $.
  Further, $\rsC $ is called \textit{connected}, if $\cC $ is a connected
  Cartan scheme.
  If $\rsC =\rsC (\cC ,(R^a)_{a\in A})$ is a root system of type $\cC $
  and $\rsC '=\rsC '(\cC ',({R'}^a_{a\in A'}))$ is a root system of
  type $\cC '$, then we say that $\rsC $ and $\rsC '$ are \textit{equivalent},
  if $\cC $ and $\cC '$ are equivalent Cartan schemes given by maps $\varphi
  _0:I\to I'$, $\varphi _1:A\to A'$ as in Def.~\ref{de:CS}, and if
  the map $\varphi _0^*:\ndZ ^I\to \ndZ ^{I'}$ given by
  $\varphi _0^*(\al _i)=\al _{\varphi _0(i)}$ satisfies
  $\varphi _0^*(R^a)={R'}^{\varphi _1(a)}$ for all $a\in A$.
\end{defin}

\begin{remar}\label{re:genrs}
  (1) Reduced root systems with a fixed basis,
  see \cite[Ch.\,VI,\S 1.5]{b-BourLie4-6},
  are examples of root systems of type $\cC $ in the following way.
  Let $\cC =\cC (I,A,(\rfl _i)_{i\in I},(\Cm ^a)_{a\in A})$
  be a Cartan scheme, such that $A=\{a\}$ has only one element, and $\Cm ^a$
  is a Cartan matrix of finite type. Then $\rfl _i=\id $ for all $i\in I$.
  Let $R^a$ be the reduced root system associated
  to $\Cm ^a$, where the basis $\{\al _i\,|\,i\in I\}$ of $\ndZ ^I$
  is identified with a basis
  of $R^a$. Then $\rsC =\rsC (\cC ,R^a)$ is a root system of type $\cC $.
  
  (2) In root systems $0$ is never a root. The same holds for the sets $R^a$
  in root systems of type $\cC $.
  Indeed, if $0\in R^a$, then $0\in R^a\cap \ndZ \al _i$ for all $i\in I$, and
  since $I$ is non-empty, this is a contradiction to (R2).

  (3) Let $\cC $ be a Cartan scheme and $\rsC $ a root system of type $\cC $.
  For $i,j\in I$ with $i\not=j$, (C1) and (R4) imply that the relations
  \begin{align}
    (\rfl _j\rfl _i)^{m_{i,j}^a}(a)=a
    \label{eq:garels}
  \end{align}
  hold for all $a\in A$. Further,
  \begin{align}
    (\s _i \s _j)^{m_{i,j}^a}1_a=(\s _j\s _i)^{m_{i,j}^a}1_a=1_a
    \label{eq:coxrel}
  \end{align}
  for all $a\in A$ and $i,j\in I$ with $i\not=j$, see Thm.~\ref{th:Coxgr}
  below. Here $1_a$ is the identity of the object $a$, and
  we use the convention that upper indices referring to objects
  are neglected if they are uniquely determined by the context.
\end{remar}

\begin{remar}\label{re:conn}
  In \cite[Def.\,2]{a-HeckYam08} it is assumed that
  a root system $\rsC $ of type $\cC $ is \textit{connected}.
  We omit this axiom to have a definition which is compatible with
  passing to restrictions, see Def.~\ref{de:rest}.
  Note that a restriction of $\rsC $ is generally not connected, even if
  $\rsC $ is connected.
\end{remar}

The following lemma states that in root systems of type $\cC $ some axioms are
redundant.

\begin{lemma}
  Let $\cC =\cC (I,A,(\rfl _i)_{i\in I},(\Cm ^a)_{a\in A})$ be a quadruple,
  where $A$ is a non-empty set, $\rho _i:A\to A$ is a map for all $i\in I$, and
  $\Cm ^a=(\cm ^a_{ij})_{i,j\in I}\in \ndZ ^{I\times I}$ for all $a\in A$,
  such that $\cm ^a_{ii}=2$ for all $i\in I$, $a\in A$, and that (C1) holds.
  For all $a\in A$ let $R^a\subset \ndZ ^I$ satisfying (R1)--(R4).
  Then $\cC $ is a Cartan scheme and
  $\rsC =\rsC (\cC ,(R^a)_{a\in A})$ is a root system of type $\cC $.
  \label{le:redund}
\end{lemma}

\begin{proof}
  We have to prove that (M1), (C2), and (M2) hold for $\Cm ^a$
  for all $a\in A$. Let $a\in A$ and $j,k\in I$ with $j\not=k$.
  Then $\al _k\in R^a$ by (R2), hence $\s _j^a(\al _k)\in R_+^{\rfl _j(a)}$
  by (R1) and (R3). Therefore $\cm ^a_{jk}\le 0$ by Eq.~\eqref{eq:sia}. This
  proves (M1) for $\Cm ^a$.

  Let now $a\in A$ and $i,j\in I$. Then
  \begin{align*}
    \s _i^{\rfl _i(a)}\s _i^a(\al _j)
    =&\s _i^{\rfl _i(a)}(\al _j-\cm ^a_{ij}\al _i)\\
    =&\al _j+(\cm ^a_{ij}-\cm ^{\rfl _i(a)}_{ij})\al _i\in R^{\rfl _i^2(a)}_+
    =R^a_+
  \end{align*}
  by Eq.~\eqref{eq:sia}, (R1)--(R3), and (C1),
  and hence $\cm ^a_{ij}\ge \cm ^{\rfl _i(a)}_{ij}$. Replacing $a$ by $\rfl
  _i(a)$, we obtain in the same way that
  $\cm ^a_{ij}\le \cm ^{\rfl _i(a)}_{ij}$. Hence (C2) holds.
  Assume now that $\cm ^a_{ij}=0$. Then $i\not=j$ by (M1). Further,
  Eq.~\eqref{eq:sia}, (C2), and relation $\cm ^a_{ij}=0$
  imply that
  \begin{align*}
    \s _i^a\s _j^{\rfl _j(a)}(\al _i)=&
    \s _i^a(\al _i-\cm ^{\rfl _j(a)}_{ji}\al _j)
    =-\al _i-\cm ^a_{ji}\al _j,
  \end{align*}
  and $\s _i^a\s _j^{\rfl _j(a)}(\al _i)\in -R^{\rfl _i(a)}_+$ by
  (R1)--(R3).
  Since $\cm ^a_{ji}\le 0$ by (M1),
  this gives that $\cm ^a_{ji}=0$, hence (M2) is proven.
\end{proof}

Recall the convention in Rem.~\ref{re:genrs}(3).
The Weyl groupoid of a root system of type $\cC $ is a generalization of the
notion of a Weyl group, as the following theorem shows.

\begin{theor}\cite[Thm.\,1]{a-HeckYam08}\label{th:Coxgr}
  Let $\cC =\cC (I,A,(\rfl _i)_{i\in I},(\Cm ^a)_{a\in A})$
  be a Cartan scheme and $\rsC =\rsC (\cC ,(R^a)_{a\in A})$ a root system
  of type $\cC $.
  Let $\Wg $ be the abstract
  groupoid with $\Ob (\Wg )=A$ such that $\Hom (\Wg )$ is
  generated by abstract morphisms $s_i^a\in \Hom (a,\rfl _i(a))$,
  where $i\in I$ and $a\in A$, satisfying the relations
  \begin{align*}
    s_i s_i 1_a=1_a,\quad (s_j s_k)^{m_{j,k}^a}1_a=1_a,
    \qquad a\in A,\,i,j,k\in I,\, j\not=k.
  \end{align*}
  Here $1_a$ is the identity of the object $a$,
  and $(s_j s_k)^\infty 1_a$ is understood to be
  $1_a$. The functor $\Wg \to \Wg (\rsC )$, which is
  the identity on the objects, and on the set of morphisms is given by
  $s _i^a\mapsto \s_i^a$ for all $i\in I$, $a\in A$, is an isomorphism of
  groupoids.
\end{theor}

One says that $\Wg (\rsC )$ is a \textit{Coxeter groupoid}.
Thus it makes sense to speak about the \textit{length}
\begin{align}
  \ell (\omega )=\min \{m\in \ndN _0\,|\,
  \omega =\s _{i_1}\cdots \s_{i_m}1_a,
  i_1,\ldots ,i_m\in I
  \}
  \label{eq:length}
\end{align}
of a morphism $\omega \in \Hom (a,b)\subset \Hom(\Wg (\rsC ))$.
The most essential difference
between Coxeter groupoids and Coxeter groups is the presence of several
objects in the former.

\begin{defin}
  Let $\cC =\cC (I,A,(\rfl _i)_{i\in I},(\Cm ^a)_{a\in A})$ be a Cartan
  scheme.
  Let $\Gamma $ be a nondirected graph,
  such that each edge is labelled by an element of $I$,
  and any two edges between two fixed vertices are labelled differently.
  Assume that there is a bijection $\varphi $ from $A$ to the set of
  vertices of $\Gamma $, and two vertices $\varphi (a)$, $\varphi (b)$, where
  $a,b\in A$, are
  connected by an edge labelled by $i\in I$ if and only if $a\not=b$ and
  $\rfl _i(a)=b$.
  The graph $\Gamma $ is called the \textit{object change diagram} of $\cC $.
  If $\rsC =\rsC (\cC ,(R^a)_{a\in A})$ is a root system of type $\cC $, then
  we also say that $\Gamma $ is the object change diagram of $\rsC $.
\end{defin}

The object change diagram of a reduced root system is a single
vertex without any edges. Other examples will appear in later sections.
Note that the object change diagram of a Cartan scheme $\cC $ is connected as
a graph if and only if the Cartan scheme $\cC $ is connected, or equivalently,
if the Weyl groupoid $\Wg (\cC )$ is a connected groupoid.

\begin{defin}
  Let $\cC =\cC (I,A,(\rfl _i)_{i\in I},(\Cm ^a)_{a\in A})$ be a Cartan scheme
  and $\rsC =\rsC (\cC ,(R^a)_{a\in A})$ a root system of type $\cC $.
  For all $a\in A$ let $(R^a)\re =\{\omega (\al _i)\,|\,\omega \in \Hom
  (b,a),b\in A, i\in I\}$, and call $(R^a)\re $ the \textit{set of real roots}
  of $a$.
\end{defin}

Note that by (R2) and (R3) we have $(R^a)\re \subset R^a$ for all $a\in A$.
The sets of real roots are interesting for various reasons, one of them is the
following.

\begin{propo}\label{pr:rsCre}
  Let $\cC =\cC (I,A,(\rfl _i)_{i\in I},(\Cm ^a)_{a\in A})$ be a Cartan
  scheme, and let $\rsC =\rsC (\cC ,(R^a)_{a\in A})$ be a root system
  of type $\cC $.
  Then $\rsC \re =\rsC \re (\cC ,((R^a)\re )_{a\in A})$
  is a root system of type $\cC $, and $\Wg (\rsC \re )=\Wg (\rsC )$.
\end{propo}

\begin{proof}
  Since $\cC $ is a Cartan scheme, it suffices to show that
  the sets $(R^a)\re $ satisfy axioms (R1)--(R4) for all $a\in A$.
  Let $a\in A$.
  Since $\omega \s _i^{\rfl _i(b)}(\al _i)=-\omega (\al _i)$ for all $i\in I$,
  $b\in A$, and $\omega \in \Hom (b,a)$, we obtain that
  $(R^a)\re =-(R^a)\re $.
  Let $(R^a)\re _+=(R^a)\re \cap \ndN _0^I$ and
  $(R^a)\re _-=(R^a)\re \cap -\ndN _0^I$. Then
  $(R^a)\re _-=-(-(R^a)\re \cap \ndN _0^I)=-(R^a)\re _+$, and hence
  (R1) implies that
  $(R^a)\re =(R^a)\re _+ \cup -(R^a)\re _+$.
  Since $\al _i\in (R^a)\re $, $(R^a)\re \subset R^a$, and $(R^a)\re
  =-(R^a)\re $ by (R1), Axiom~(R2) implies that
  $(R^a)\re \cap \ndZ \al _i=\{\al _i,-\al _i\}$.
  Axiom~(R3) holds for $\rsC \re $ by definition. Finally, if $m^a_{i,j}$
  is finite for an $a\in A$ and elements $i,j\in I$ with $i\not=j$,
  then
  \[R^a\cap (\ndN _0\al _i+\ndN _0\al _j)=(R^a)\re \cap (\ndN _0\al
  _i+\ndN _0\al _j)\]
  by \cite[Lemma\,4]{a-HeckYam08}. Thus
  (R4) holds for $\rsC \re $, since it holds for $\rsC $. The equation
  $\Wg (\rsC \re )=\Wg (\rsC )$ follows since $\Wg (\rsC \re )=\Wg (\cC )=\Wg
  (\rsC )$ by definition.
\end{proof}

Now we discuss the finiteness of root systems of type $\cC $.

\begin{defin}
  Let $\cC =\cC (I,A,(\rfl _i)_{i\in I},(\Cm ^a)_{a\in A})$ be a Cartan scheme
  and $\rsC =\rsC (\cC ,(R^a)_{a\in A})$ a root system of type $\cC $.
  We say that $\rsC $ is \textit{finite} if $R^a$ is finite for all $a\in A$.
\end{defin}

The finiteness of $\rsC $ does not mean that $\Wg (\rsC )$ is finite, since
$A$ may be infinite. But the following holds.

\begin{lemma}
  Let $\cC =\cC (I,A,(\rfl _i)_{i\in I},(\Cm ^a)_{a\in A})$ be a
  connected Cartan scheme
  and $\rsC =\rsC (\cC ,(R^a)_{a\in A})$ a root system of type $\cC $.
  Then the following are equivalent.
  \begin{enumerate}
    \item $\rsC $ is finite.
    \item $R^a$ is finite for at least one $a\in A$.
    \item $\rsC \re $ is finite.
    \item $\Wg (\rsC )$ is finite.
  \end{enumerate}
  \label{le:finHom}
\end{lemma}

\begin{proof}
  (1)$\Rightarrow $(2) is trivial.
  
  (2)$\Rightarrow $(1). Assume that $a\in A$ such that $R^a$ is finite.
  Since $\rsC $ is connected, for each $b\in A$ there exists $\omega \in
  \Hom (a,b)$. Then $R^b=\omega (R^a)$ by (R3), and hence $R^b$ is finite for
  all $b\in A$.
  
  (1)$\Rightarrow $(4). Let $a,b\in A$.
  Since $R^b$ is finite, $R^a$ contains the standard basis of $\ndZ ^I$ by
  (R2), and since $\omega (R^a)\subset R^b$
  for all $\omega \in \Hom (a,b)$, the set $\Hom (a,b)$ is
  finite. Assume now that $A$ is infinite, and let $a\in A$. The finiteness
  of $R^a$ implies that there exist $b,c\in A$ and $f\in \Hom (a,b)$,
  $g\in \Hom (a,c)$ with $f\not=g$ and
  \begin{align*}
    Y:=\{\al \in R_+^a\,|\,f(\al )\in R_+^b\}=
    \{\beta \in R_+^a\,|\,g(\beta )\in R_+^c\}.
  \end{align*}
  Then $fg^{-1}\in \Hom (c,b)$ and
  \begin{align*}
    fg^{-1}(g(Y))=&f(Y)\subset R_+^b,\\
    fg^{-1}(R_+^c\setminus g(Y))=&f(-(R_+^a\setminus Y))\subset R_+^b
  \end{align*}
  by (R1) and (R3).
  Therefore $fg^{-1}(R^c_+)\subset R^b_+$, and hence
  $b=c$ and $fg^{-1}=1_b$ by \cite[Lemma~8(iii)]{a-HeckYam08}. This is a
  contradiction to $f\not=g$, and hence $A$ is finite. This proves
  (1)$\Rightarrow $(4).

  (4)$\Rightarrow $(1). We prove that if $\rsC $ is infinite, then $\Wg
  (\rsC )$ is infinite. For this we show by induction on $m$ that
  for all $m\in \ndN $ there exist $a,b\in A$ and $\omega \in \Hom
  (a,b)$ such that
  \begin{align}\label{eq:posneg}
    |\{\al \in R_+^a\,|\,\omega (\al )\in - R_+^b\}|=m.
  \end{align}
  The latter holds for $m=1$, since $\omega =\s _i^a$ fulfills
  Eq.~\eqref{eq:posneg} for all $a\in A$, $i\in I$, and $b=\rfl _i(a)$.
  Suppose now that $m\in \ndN $, $a,b\in A$, and $\omega \in \Hom (a,b)$
  such that Eq.~\eqref{eq:posneg} holds. Since $|R^a|=\infty $, there exists
  $\al \in R_+^a$ with $\omega (\al )\in R_+^b$.
  Since $\omega $ is linear, there exists
  $i\in I$ such that $\omega (\al _i)\in R_+^b$.
  Then $\ell (\omega \s _i^{\rfl _i(a)})=m+1$
  by \cite[Cor.\,3]{a-HeckYam08}, and hence
  \[\big|\{\al \in R_+^{\rfl _i(a)}\,|\,\omega \s _i^{\rfl _i(a)}(\al )\in
  - R_+^b\}\big|=m+1\]
  by \cite[Lemma\,8(iii)]{a-HeckYam08}. Thus the induction step is proved.

  Finally, the equivalence of (3) and (4) follows from the equivalence of (1)
  and (4) and the equation $\Wg (\rsC \re )=\Wg (\rsC )$, see
  Prop.~\ref{pr:rsCre}.
\end{proof}

Now we prove that if $\rsC $ is a finite root system of type $\cC $,
then all roots are real, that is,
$\rsC $ is uniquely determined by $\cC $.

\begin{propo} \label{pr:allposroots}
  Let $\cC =\cC (I,A,(\rfl _i)_{i\in I},(\Cm ^a)_{a\in A})$ be a Cartan scheme
  and $\rsC =\rsC (\cC ,(R^a)_{a\in A})$ a root system of type $\cC $.
  Let $a\in A$, $m\in \ndN _0$, and $i_1,\ldots ,i_m\in I$ such that
  $\om =1_a\s _{i_1}\s _{i_2}\cdots \s_{i_m}$ and $\ell (\om )=m$.
  Then the elements
  \begin{align*}
    \beta _n=1_a\s _{i_1}\s _{i_2}\cdots \s _{i_{n-1}}(\al _{i_n})\in R^a_+,
  \end{align*}
  where $n\in \{1,2,\ldots ,m\}$ (and $\beta _1=\al _{i_1}$),
  are pairwise different. In particular, if $\rsC $ is finite and
  $\om \in \Hom (\Wg (\rsC ))$ is a longest element,
  see \cite[Cor.\,5]{a-HeckYam08}, then
  \begin{align*}
    \{\beta _n\,|\,1\le n\le \ell (\om )=|R^a|/2\}=R_+^a.
  \end{align*}
\end{propo}

\begin{proof}
  For all $n\in \ndN _0$ with $n\le m$ let
  $a_n = \rfl _{i_n}\cdots \rfl _{i_2}\rfl _{i_1}(a)$.
  Since $\ell (\s _{i_1}^{a_1}\s _{i_2}^{a_2}\cdots \s _{i_n}^{a_n})=n$
  for all $n\in \{1,2,\ldots ,m\}$,
  one has $\beta _n\in R_+^a$ for these $n$ by \cite[Cor.\,3]{a-HeckYam08}.
  By the same argument one has for all $k,n\in \ndN $ with $k<n\le m$ the
  relation
  \begin{align*}
    \s _{i_{k+1}}^{a_{k+1}}\s _{i_{k+2}}^{a_{k+2}}\cdots
    \s _{i_{n-1}}^{a_{n-1}}(\al _{i_n})
    \in R_+^{a_k}.
  \end{align*}
  Thus $\s _{i_k}\s _{i_{k+1}}\cdots \s _{i_{n-1}} 1_{a_{n-1}}(\al _{i_n})\not=
  \al _{i_k}$, and hence $\beta _n\not=\beta _k$ for all $k,n\in \ndN $ with
  $k<n\le m$.
\end{proof}

For any groupoid $G$ and any $a\in \Ob (G)$ let
$\Hom (a)=\Hom (a,a)\subset \Hom (G)$.
Then $\Hom (a)$ is a subgroup of $G$, which depends on $a$. However, the
following is true.

\begin{propo}
  Let $G$ be a connected groupoid and $a,b\in \Ob (G)$.
  Then $\Hom (a)$ and $\Hom (b)$ are isomorphic groups.
\end{propo}

\begin{proof}
  Choose $X_b\in \Hom (a,b)$. This exists since $G$ is connected.
  Then the map
  \begin{align*}
    \phi _{a,b}:\Hom (a)\to \Hom (b),\quad
    g\mapsto X_b g X_b^{-1}
  \end{align*}
  is a group homomorphism with inverse given by
  $\phi _{a,b}^{-1}(g)=X_b^{-1}g X_b$.
\end{proof}

The map $\phi _{a,b}$ in the previous proof is a piece of a more general
structure. Namely, let $G$ be a connected groupoid, $a\in \Ob (G)$, and
for each $b\in \Ob (G)$ let $X_b\in \Hom (a,b)$ be a fixed morphism.
Then the assignment $F_{a,X}:G\to \Hom (a)$,
\begin{equation}
\begin{aligned}
  F_{a,X}(b)=&a & &\text{for all $b\in \Ob (G)$,}\\
  F_{a,X}(g)=&X_c^{-1}gX_b & &\text{for all $g\in \Hom (b,c)$,}
\end{aligned}
  \label{eq:grpidfunctor}
\end{equation}
defines a fully faithful functor. In fact, $G$ is as a groupoid isomorphic to
the transformation groupoid $H:=\Stab (a)\times \Ob (G)$ given by $\Ob (H)=\Ob
(G)$, $\Hom (H)=\Ob (G)\times \Stab (a)\times \Ob (G)$ with composition
\begin{align*}
  (b,g,c)(b',g',c')=
  \begin{cases}
    \text{not defined} & \text{if $c\not=b'$,}\\
    (b,gg',c') & \text{if $c=b'$.}
  \end{cases}
\end{align*}
The isomorphism $G\to H$ is given by $g\mapsto (c,F_{a,X}(g),b)$
for $g\in \Hom (b,c)$. In particular, $G$ is uniquely determined by
the cardinality of $\Ob (G)$ and by $\Stab (a)$ for any $a\in \Ob (G)$.

If a connected groupoid $G$ is presented by generators and relations, then
for any $a\in \Ob (G)$ the group $\Hom (a)$ also can be presented by
generators and relations. To do so, let $F_{a,X}:G\to \Stab (a)$ be the
functor defined above. The following proposition then follows from the
discussion above.

\begin{propo}\label{pr:Stabrel}
  Let $G$ be a connected groupoid and let $a\in \Ob (G)$. Suppose that
  $J,K$ are index sets and $\Hom (G)$ is generated by $s_j\in \Hom (a_j,b_j)$,
  where $j\in J$, and relations $r_k=1_{c_k}\in \Hom (c_k)$, where $k\in K$.
  Then $\Hom (a)$ is generated by $F_{a,X}(s_j)$, where $j\in J$,
  and relations $F_{a,X}(r_k)=1$, where $k\in K$ and $1$ is the neutral
  element of $\Hom (a)$.
\end{propo}

\section{Standard Cartan schemes and their Weyl groupoids}

In this section we study root systems of type $\cC $,
where the Cartan matrices are identical for all objects.
The structure of more general root systems of type $\cC $ seems to be much
more complicated, as the classification results in the next sections show.

\begin{defin}\label{de:standard}
  Let $\cC =\cC (I,A,(\rfl _i)_{i\in I},(\Cm ^a)_{a\in A})$ be a Cartan
  scheme such that $\Cm ^a=\Cm ^b$ for all $a,b\in A$. Then we say that $\cC $
  is \textit{standard}, and that the Weyl groupoid $\Wg (\cC )$ is standard.
  If $\rsC $ is a root system of type $\cC $, then we say that $\rsC $ is
  standard, if $\cC $ is a standard Cartan scheme.
\end{defin}

The standard Cartan schemes in the next example show that the class of root
systems of type $\cC $, where $\cC $ is running over all Cartan schemes,
is richer than the one of finite groups.

\begin{examp}\label{ex:rcwithgivenHom}
  Let $H$ be a finite group. Then there exists a connected Cartan scheme
  $\cC =\cC (I,A,(\rfl _i)_{i\in I},(\Cm ^a)_{a\in A})$ and a finite root
  system $\rsC =\rsC (\cC ,(R^a)_{a\in A})$ of type $\cC $,
  such that $\Stab (a)\cong H$ for all $a\in A$.
  Indeed, $H$ can be considered as a subgroup of a symmetric group $\SG{n+1}$.
  Let $A$ be the set of left cosets $gH$, where $g\in \SG{n+1}$, and let
  $I=\{1,2,\ldots ,n\}$.
  For all $gH\in A$ and $i\in I$ let $\rfl _i(a)=(i,i+1)gH$,
  where $(i,i+1)$ is the transposition of $i$ and $i+1$, and
  $\Cm ^{gH}=(\cm ^{gH}_{i j})_{i,j\in I}$ with
  \begin{align*}
    \cm ^{gH}_{ij}=
    \begin{cases}
      2 & \text{if $i=j$,}\\
      -1 & \text{if $|i-j|=1$,}\\
      0 & \text{otherwise.}
    \end{cases}
  \end{align*}
  Then $\cC =\cC (I,A,(\rfl _i)_{i\in I},(\Cm ^a)_{a\in A})$ is a standard
  Cartan scheme. For all $a\in A$ let
  $R^a=R_+^a\cup - R_+^a$, where
  \begin{align}
    R_+^a:=\{\al _i+\al _{i+1}+\cdots +\al _j\,|\,
    1\le i\le j\le n\},
    \label{eq:Snroots}
  \end{align}
  be the set of roots associated to $\SG{n+1}$. Then
  $\rsC =\rsC (\cC ,(R^a)_{a\in A})$ is a root system of type $\cC $,
  and $\Hom (eH)\subset \Hom (\Wg (\rsC ))$ is isomorphic to $H$.
\end{examp}

The structure of finite connected standard root systems of type $\cC $ is
very close to the structure of reduced root systems in the sense of
\cite[Ch.\,VI,\,\S 1.4]{b-BourLie4-6}.

\begin{theor}\label{th:standard}
  Let $I$ be a non-empty finite set, $\Cm =(\cm _{i j})_{i,j\in I}$ a
  generalized Cartan matrix, and
  $\cC =\cC (I,A,(\rfl _i)_{i\in I},(\Cm ^a)_{a\in A})$
  a connected standard Cartan scheme with $\Cm ^a=\Cm $ for all $a\in A$,
  and let $\rsC =\rsC (\cC ,(R^a)_{a\in A})$ be
  a root system of type $\cC $.

  (1) For all $a\in A$ the set $\cup _{b\in A}\Hom (a,b)\subset \Aut (\ndZ
  ^\theta )$ is a group, and as such it is isomorphic
  to the Weyl group $W(\Cm )$ associated to the generalized Cartan matrix
  $\Cm $.
  
  (2) $\rsC $ is finite if and only if $\Cm $ is of finite type.

  (3) Assume that $\rsC $ is finite. Then for all $a\in A$,
  $R^a$ is the set of roots corresponding to $W(\Cm )$,
  and hence independent of the choice of $a\in A$.
\end{theor}

\begin{proof}
  Since $\cC $ is standard,
  the maps $\s _i^a\in \Aut (\ndZ ^\theta )$ do not depend on the object $a\in
  A$, and generate the Weyl group $W(\Cm )\subset \Aut (\ndZ ^\theta )$
  associated to the generalized Cartan matrix $\Cm $. Let $a\in A$.
  Since
  \begin{align*}
    \cup _{b\in A}\Hom (a,b)=\{(\rfl _{i_1}\cdots \rfl _{i_n}(a),
    &\s _{i_1}\cdots \s _{i_n}1_a,a)\,|\\
    &n\in \ndN _0,\,i_1,\ldots ,i_n\in I\},
  \end{align*}
  Thm.~\ref{th:Coxgr} implies that (1) holds.

  Assume that $\rsC $ is finite.
  Since $\cC $ is standard, Prop.~\ref{pr:allposroots} tells that
  $R^a_+$ is the set of positive roots corresponding to $W(\Cm )$.
  This implies (3).

  Now we prove (2). If $\Cm $ is of finite type, then $W(\Cm )$ is finite.
  Since $\cC $ is connected, $A$ is finite by Part~(1). Thus $\Wg (\cC )$ is
  finite by Part~(1), and hence $\rsC $ is finite by Lemma~\ref{le:finHom}.

  Conversely, assume that $\rsC $ is finite. Then Lemma~\ref{le:finHom}
  implies that $\Wg (\rsC )$ is finite,
  and hence $W(\Cm )$ is finite by Part~(1). Thus
  $\Cm $ is of finite type.
\end{proof}

\section{Decomposition of finite root systems}
\label{sec:dec}

In this section we study the reducibility of root systems of type $\cC $.
An analogous notion exists for root systems, see
\cite[Ch.\,6,\,\S 1.2]{b-BourLie4-6},
and it is crucial for classification results.

\begin{defin}\label{de:rest}
Let $\cC =\cC (I,A,(\rfl _i)_{i\in I}, (\Cm ^a)_{a\in A})$
be a Cartan scheme.
Let $J\subset I$ be a non-empty subset, and identify
$\{\al _i\,|\,i\in J\}$ with the standard basis of $\ndZ ^J$.
For all $a\in A$ let
${\Cm '}^a=(\cm ^a_{ij})_{i,j\in J}$.
Then $\cC '=\cC '(J,A,(\rfl _i)_{i\in J},({\Cm '}^a)_{a\in A})$ is a Cartan
scheme, called the \textit{restriction of} $\cC $ to $J$, and will be denoted
by $\cC |_J$.

Let $\rsC =\rsC (\cC ,(R^a)_{a\in A})$ be
a root system of type $\cC $. Define
${R'}^a=R^a\cap \sum _{i\in J}\ndZ \al _i$. Then
$\rsC '=\rsC '(\cC |_J,({R'}^a)_{a\in A})$ is a root system of type $\cC |_J$,
and will be denoted by $\rsC |_J$.
\end{defin}

Restrictions are helpful to decide if a root system of type $\cC $ is standard.

\begin{remar}\label{re:standard}
  Let $\cC =\cC (I,A,(\rfl _i)_{i\in I}, (\Cm ^a)_{a\in A})$
  be a Cartan scheme and $\rsC =\rsC (\cC ,(R^a)_{a\in A})$
  a root system of type $\cC $.
  Assume that for each pair $(i,j)\in I\times I$ with
  $i\not=j$ there exists a subset $J\subset I$ such that $i,j\in J$ and $\rsC
  |_J$ is standard. Then $\rsC $ is standard.
  Indeed, $\rsC $ is standard if and only if
  $\cm ^a_{ij}=\cm ^b_{ij}$ for all $a,b\in A$ and $i,j\in I$.
  The latter holds by assumption on the pairs $(i,j)\in I\times I$,
  and since $\cm ^a_{ii}=2$ for all $i\in I$ and $a\in A$.
\end{remar}

\begin{defin}\label{de:dirsum}
  Let $\cC =\cC (I,A,(\rfl _i)_{i\in I},({\Cm }^a)_{a\in A})$ be a  Cartan
  scheme. Assume that $I',I''\subset I$ are non-empty disjoint subsets such
  that $I=I'\cup I''$ and $\cm ^a_{i j}=0$ for all $i\in I'$, $j\in I''$.
  Then we write $\cC =\cC |_{I'}\oplus \cC |_{I''}$, and say that $\cC $ is
  the {\em direct sum of} $\cC |_{I'}$ {\em and} $\cC |_{I''}$.
  
  Let $\rsC =\rsC (\cC ,(R^a)_{a\in A})$ be a root system of type $\cC $.
  Assume that
  \[ R^a=\Big(R^a\cap \sum _{i\in I'}\ndZ \al _i\Big)\cup
  \Big(R^a\cap \sum _{j\in I''}\ndZ \al _j\Big)\qquad
  \text{for all $a\in A$.}\]
  Then we write $\rsC =\rsC |_{I'}\oplus \rsC |_{I''}$, and $\rsC $ is called the
  {\em direct sum of} $\rsC |_{I'}$ {\em and} $\rsC |_{I''}$.
  We also say that $\rsC $ is {\em reducible}.
  If $\rsC \not=\rsC |_{I_1}\oplus \rsC |_{I_2}$ for all non-empty disjoint
  subsets $I_1,I_2\subset I$, then $\rsC $ is termed {\em irreducible}.
\end{defin}

{}From now on let $\cC =\cC (I,A,(\rfl _i)_{i\in I}, (\Cm ^a)_{a\in A})$
be a connected Cartan scheme and $\rsC =\rsC (\cC ,(R^a)_{a\in A})$
a root system of type $\cC $.
We are going to give criteria for the reducibility of $\rsC $.

\begin{lemma}
  Let $a\in A$ and $i,j\in
  I$ with $i\not=j$. The following are equivalent.
  \begin{enumerate}
    \item $\cm ^a_{ij}=\cm ^a_{ji}=0$.
    \item $R^a\cap (\ndN _0\al _i+\ndN _0\al _j)=\{\al _i,\al _j\}$.
    \item $m_{i,j}^a=2$.
  \end{enumerate}
  \label{le:m=2}
\end{lemma}

\begin{proof}
  Since $\al _i,\al _j\in R^a$ by (R2), (2) is equivalent to (3).
  Further, from (R1)--(R3) and Eq.~\eqref{eq:sia} we conclude that (2) implies
  (1). Assume now that (1) holds, and let $\al :=r_i\al _i+r_j\al _j\in R^a_+$,
  where $r_i,r_j\in \ndN _0$. Then (R1) and relation
  $\s _i^a(\al )=-r_i\al _i+r_j\al _j\in R^{\rfl _i(a)}$
  imply that $r_i=0$ or $r_j=0$.
  Hence $\al \in \{\al _i,\al _j\}$ by (R2). This proves (1)$\Rightarrow $(2).
\end{proof}

\begin{lemma}
  Suppose that $a\in A$ and $i,j\in I$, where $i\not=j$.
  If $\cm ^a_{i j}=0$ then $\cm ^a_{j l}=\cm ^{\rfl _i(a)}_{j l}$
  for all $l\in I$.
  \label{le:samerowfor0}
\end{lemma}

\begin{proof}
  Let $l\in I$.
  If $l=j$, then $\cm ^a_{j l}=\cm ^{\rfl _i(a)}_{j l}=2$, and if $l=i$, then
  $\cm ^a_{j l}=\cm ^a_{l j}=0=\cm ^{\rfl _i(a)}_{l j}=\cm ^{\rfl _i(a)}_{j l}$
  by (M2) and (C2). Assume now that $l\in I\setminus \{i,j\}$. Then
  $\s _i^{\rfl _j(a)}\s _j^a(\al _l)=\s _j^{\rfl _i(a)}\s _i^a(\al _l)$
  by Thm.~\ref{th:Coxgr}. Explicit calculation gives
  \begin{align*}
    \s _i^{\rfl _j(a)}\s _j^a(\al _l)
    =\s _i^{\rfl _j(a)}(\al _l-\cm ^a_{j l}\al _j)
    =\al _l-\cm ^{\rfl _j(a)}_{i l}\al _i-\cm ^a_{j l}\al _j,
  \end{align*}
  and similarly
  $\s _j^{\rfl _i(a)}\s _i^a(\al _l)
   =\al _l-\cm ^{\rfl _i(a)}_{j l}\al _j-\cm ^a_{i l}\al _i$. Comparing
   the coefficients of $\al _j$ gives the claim of the lemma.
\end{proof}

\begin{propo}\label{pr:decomp}
  Let $I'\subset I$ be a subset of $I$, and let $I''=I\setminus I'$.
  Assume that $I',I''\not=\emptyset $.
  The following are equivalent.
\begin{enumerate}
  \item There exists $a\in A$ such that $\cm ^a_{i j}=0$ for all $i\in I'$ and
    $j\in I''$.
  \item For all $a\in A$ and $i\in I'$, $j\in I''$ one has
    $\cm ^a_{i j}=\cm ^a_{j i}=0$.
  \item
    For all $a\in A$ let $\hat{R}^a=(R^a\cap \sum _{i\in I'}\ndZ \al _{i'})
    \cup (R^a \cap \sum _{i\in I''}\ndZ \al _{i''})$.
    Then $\hat{\rsC }=\hat{\rsC }(\cC ,(\hat{R}^a)_{a\in A})$
    is a root system of type $\cC $.
\end{enumerate}
  If $\rsC $ is finite then (1)--(3) are equivalent to the following.
  \begin{enumerate}
    \item[(4)] $\hat{\rsC }=\rsC $, where $\hat{\rsC }$ is as in \textrm{(3)}.
    \item[(5)] $\rsC =\rsC |_{I'}\oplus \rsC |_{I''}$ with respect to the
      permutation $\phi =\id $ of $A$.
  \end{enumerate}
\end{propo}

\begin{proof}
  The implication (2)$\Rightarrow $(1) is trivial.

  (1)$\Rightarrow $(2). Let $b\in A$ and $l\in I'$, and assume that
  $\cm ^b_{ij}=\cm ^b_{ji}=0$ for all $i\in I'$, $j\in I''$.
  Since $\cm ^b_{lj}=0$ for all $j\in I''$,
  Lemma~\ref{le:samerowfor0} implies that $\cm ^{\rfl _l(b)}_{j i}=\cm ^b_{j
  i}=0$ for all $i\in I'$ and $j\in I''$, that is, $\cm ^{\rfl _l(b)}_{i j}
  =\cm ^{\rfl _l(b)}_{j i}=0$ for all $i\in I'$, $j\in I''$. Together with the
  analogous argument for $l\in I''$ we obtain that
  $\cm ^{\rfl _l(b)}_{i j}=\cm ^{\rfl _l(b)}_{j i}=0$ for all $i\in I'$,
  $j\in I''$, and $l\in I$. Thus (2) follows from (1) since $\rsC $
  is connected.

  (3)$\Rightarrow $(2). Since $\hat{R}^a\cap (\ndN _0\al _i+\ndN _0\al _j)=
  \{\al _i,\al _j\}$ for all $i\in I'$, $j\in I''$, and $a\in A$,
  (2) follows from Lemma~\ref{le:m=2} and (3).

  (2)$\Rightarrow $(3).
  Since $\rsC $ is a root system of type $\cC $,
  (2) and Lemma~\ref{le:m=2} imply that $m^a_{i,j}=2$ for all $a\in A$, $i\in
  I'$, and $j\in I''$. Thus (3) is equivalent to the
  fact that $\s _l^a(\hat{R}^a)\subset \hat{R}^{\rfl _l(a)}$
  for all $l\in I$ and $a\in A$.
  Let $\al \in \hat{R}^a$. Then (2) implies that
  \begin{align*}
    \s _l^a(\al )\in R^{\rfl _l(a)}\cap \big(\sum _{i'\in I'}\ndZ \al _{i'}
    \cup \sum _{i''\in I''}\ndZ \al _{i''}\big)
    =\hat{R}^{\rfl _l(a)},
  \end{align*}
  and hence (3) holds.

  Assume now that $\rsC $ is finite.
  Then, by Prop.~\ref{pr:allposroots}, (2) implies that
  $\hat{R}^a_+=R^a_+$ for all $a\in A$, that is, (4) holds.
  Obviously, (4) implies (3), and (4) is also equivalent to (5),
  hence the proposition is proven.
\end{proof}

For the equivalence $(3)\Leftrightarrow (4)$ in Prop.~\ref{pr:decomp}
the finiteness assumption on $\rsC $ is necessary, as
the following example shows.

\begin{examp}
  Let $I=\{1,2,3,4\}$, $A=\{a\}$, and
  \begin{gather*}
    \Cm ^a=
    \begin{pmatrix}
      2 & -2 & 0 & 0 \\
      -2 & 2 & 0 & 0 \\
      0 & 0 & 2 & -2 \\
      0 & 0 & -2 & 2
    \end{pmatrix},\\
    \begin{aligned}
    R^a_+=\{&
    m\al _1+(m+1)\al _2,
    (m+1)\al _1+m\al _2,\\
    &m\al _3+(m+1)\al _4,
    (m+1)\al _3+m\al _4,
    \,|\,m\in \ndN _0\}\\
    \cup \{&\al _1+\al _2+\al _3+\al _4\}.
    \end{aligned}
  \end{gather*}
  Then (2) holds in Prop.~\ref{pr:decomp} for $I'=\{1,2\}$ and $I''=\{3,4\}$.
  However, $R^a\not=R^a|_{\{1,2\}}\cup R^a|_{\{3,4\}}$ since
  $\al _1+\al _2+\al _3+\al _4\notin \hat{R}^a$, and hence $\rsC $ is
  irreducible.
\end{examp}

We continue with some general statements about $\rsC $.

\begin{lemma}
  Suppose that $a\in A$ and $i,j\in I$, where $i\not=j$.
  The following are equivalent.
  \begin{enumerate}
    \item $\cm ^a_{i j}=\cm ^a_{j i}=-1$.
    \item $R^a\cap (\ndN _0\al _i+\ndN _0\al _j)=
      \{\al _i, \al _i+\al _j,\al _j\}$.
    \item $m_{i,j}^a=3$.
  \end{enumerate}
  \label{le:Cartan-1}
\end{lemma}

\begin{proof}
  To (1)$\Rightarrow $(2).
  Let $\al =c_1\al _i+c_2\al _j$ with $c_1,c_2\in \ndN $, and assume that
  $\al \in R^a$.
  Then $\s _i^a(\al )=(c_2-c_1)\al _i+c_2\al _j$,
  and hence relation $\s _i^a(\al )\in R^{\rfl _i(a)}_+\cup - R^{\rfl _i(a)}_+$
  tells that $c_2\ge c_1$. By symmetry one gets $c_1=c_2$,
  and hence $\s _i^a(\al )=c_2\al _j$.
  By (R2) one obtains that $c_2=1$.

  The implication (2)$\Rightarrow $(3) follows from the definition of
  $m_{i,j}^a$. We have to prove (3)$\Rightarrow $(1).
  Assume that $m_{i,j}^a=3$. Then $\cm ^a_{i j}$, $\cm ^a_{ji}<0$ by
  Lemma~\ref{le:m=2}. Thus,
  $\s _i^{\rfl _i(a)}(\al _j)=\al _j-\cm ^{\rfl _i(a)}_{i j}\al _i
  \in R_+^a\setminus \{\al _i,\al _j\}$, and hence
  $\beta _1:=\al _j-\cm ^a_{ij}\al _i\in R_+^a\setminus \{\al _i,\al _j\}$
  by (C2).
  Similarly, $\beta _2:=\al _i-\cm ^a_{ji}\al
  _j\in R_+^a\setminus \{\al _i,\al _j\}$, and therefore (3) implies that
  $\beta _1=\beta _2$, that is, (1) holds.
\end{proof}

\begin{lemma}
  Suppose that $a\in A$ and $i,j\in I$ such that $m_{i,j}^a=3$.
  Then $\cm^{\rfl _i(a)}_{il}+\cm ^{\rfl _i(a)}_{jl}
  =\cm ^{\rfl _i\rfl _j(a)}_{il}+\cm ^{\rfl _i\rfl _j(a)}_{jl}$
  for all $l\in I$.
  \label{le:m=3}
\end{lemma}

\begin{proof}
  If $l\in \{i,j\}$, then Lemma~\ref{le:Cartan-1} implies that 
  both sides of the claimed equation are equal to $1$. Assume now that
  $l\in I\setminus \{i,j\}$. Then
  \[\s _i^{\rfl _j\rfl _i(a)}\s _j^{\rfl _i(a)}\s _i^a(\al _l)=
  \s _j^{\rfl _i\rfl _j(a)}\s _i^{\rfl _j(a)}\s _j^a(\al _l)\]
  by Thm.~\ref{th:Coxgr}, that is,
  \begin{align*}
  \al _l-(\cm ^{\rfl _i(a)}_{jl}+\cm ^{\rfl _j\rfl _i(a)}_{il})\al _i
  &-(\cm ^{\rfl _i(a)}_{jl}+\cm ^a_{il})\al _j\\
  =&
  \al _l
  -(\cm ^{\rfl _j(a)}_{il}+\cm ^a_{jl})\al _i
  -(\cm ^{\rfl _j(a)}_{il}+\cm ^{\rfl _i\rfl _j(a)}_{jl})\al _j.
  \end{align*}
  One obtains the claim of the lemma by comparing the coefficients of $\al _j$
  and by using (C2).
\end{proof}

\begin{lemma}
  \label{le:Cartan3}
  Suppose that $\rsC $ is finite. Let $a,b\in A$
  and $i,j,l\in I$ such that $i\not=j$,
  $\rfl _i(a)=\rfl _j(a)=b\not=a$, and $\rfl _l(a) =a$.
  Then the following hold.
  \begin{enumerate}
    \item $\cm ^a_{i j}\cm ^a_{i l}\cm ^a_{j l}=0$.
    \item If $\rfl _i\rfl _l(b)\not=\rfl _l(b)$ and
      $\rfl _j\rfl _l(b)\not=\rfl _l(b)$ then
      $\cm ^a_{l n}=\cm ^b_{l n}$ for all $n\in I$.
  \end{enumerate}
\end{lemma}

\begin{proof}
  Note that $\cm ^b_{j i}=\cm ^a_{j i}$ and $\cm ^b_{j l}=\cm ^a_{j l}$ by
  (C2).

  Let $\tilde{\s }=\s _l^a\s _j^b\s _i^a\in \Hom (a)$. Since $\rsC $ is finite,
  $\tilde{\s }$ must have finite order by Lemma~\ref{le:finHom}.
  Let $W=\Lin _\ndZ \{\al _i,\al _j,\al _l\}$. Since $\tilde{\s }(W)\subset W$,
  $\s :=\tilde{\s }|_W \in \End (W)$ has to have finite order as well.
  This will yield both Claim~(1) and Claim~(2).

  To (1). One has $\s (\al _i)\in -\al _i
  -\ndN _0\al _j-\ndN _0\al _l$. If $\s (\al _i)=-\al _i$, then
  $\cm ^a_{l i}=\cm ^a_{i l}=0$ and hence (1) holds.
  
  It remains to consider the case when $\s \not=\pm \id $.
  Then finiteness of the order of $\s $ tells that
  $|\tr \s |\in \{0,1,2\}$. Explicit calculation gives that
  \begin{align}
    \tr \s = -\cm ^a_{i l}\cm ^a_{j i}\cm ^a_{l j}+\cm ^a_{i j} \cm ^a_{j i}
    +\cm ^a_{i l}\cm ^a_{l i}+\cm ^a_{j l}\cm ^a_{l j}-3.
    \label{eq:trs}
  \end{align}
  Exchanging $i$ and $j$ further implies that
  \begin{align*}
    -2\le -\cm ^a_{i j}\cm ^a_{j l}\cm ^a_{l i}+\cm ^a_{i j} \cm ^a_{j i}
    +\cm ^a_{i l}\cm ^a_{l i}+\cm ^a_{j l}\cm ^a_{l j}-3\le 2.
  \end{align*}
  If $\cm ^a_{i j}\cm ^a_{i l}\cm ^a_{j l}\not=0$,
  then the latter relations and (M1), (M2) imply that
  \begin{align*}
    \cm ^a_{i j}=\cm ^a_{j i}=
    \cm ^a_{i l}=\cm ^a_{l i}=
    \cm ^a_{j l}=\cm ^a_{l j}=-1.
  \end{align*}
  In this case one has
  \begin{align*}
    \s =
    \begin{pmatrix}
      -1 & 1 & 1 \\ -1 & 0 & 2 \\ -2 & 1 & 2
    \end{pmatrix},\quad
    \s ^2=
    \begin{pmatrix}
      -2 & 0 & 3 \\ -3 & 1 & 3 \\ -3 & 0 & 4
    \end{pmatrix}
  \end{align*}
  with respect to the basis $\{\al _i,\al _j,\al _l\}$.
  Then $\tr \s ^2=3$, but $\s ^2\not=\id $,
  and hence $\s $ does not have finite order. This is a contradiction to the
  assumption that $\rsC $ is finite, and hence (1) holds.

  To (2). Assume that $\rfl _i\rfl _l(b)\not=\rfl _l(b)$ and
  $\rfl _j\rfl _l(b)\not=\rfl _l(b)$.
  If $\cm ^a_{i l}\cm ^a_{j l}=0$ then Claim (2) is valid by
  Lemma~\ref{le:samerowfor0}. Thus by (1) it remains to consider the setting
  $\cm ^a_{i j}=0$, $\cm ^a_{i l}\cm ^a_{j l}\not=0$. In this case one has
  \begin{align*}
    \s =
    \begin{pmatrix}
      -1 & 0 & -\cm ^a_{i l} \\ 0 & -1 & -\cm ^a_{j l} \\
      \cm ^a_{l i} & \cm ^a_{l j} &
      t-1
    \end{pmatrix},\quad
    \s ^2=
    \begin{pmatrix}
      1-\cm ^a_{i l}\cm ^a_{l i} & -\cm ^a_{i l}\cm ^a_{l j} & * \\
      * & 1-\cm ^a_{j l}\cm ^a_{l j} & * \\
      * & * & t^2-3t+1
    \end{pmatrix}
  \end{align*}
  with respect to the basis $\{\al _i,\al _j,\al _l\}$,
  where $t=\cm ^a_{i l}\cm ^a_{l i}+\cm ^a_{j l}\cm ^a_{l j}$.
  In the next paragraph we prove that $t\ge 4$. This implies that
  \begin{align*}
    \tr \s ^2=t^2-4t+3=(t-2)^2-1\ge 3.
  \end{align*}
  Since $\s ^2\not=\id $, part (2) of the lemma is proven.

  Recall that $\cm ^a_{i l}\not=0$. Further,
  relation $\rfl _i\rfl _l(b)\not=\rfl _l(b)$ implies that
  $\rfl _i\rfl _l\rfl _i(a)\not=\rfl _l\rfl _i\rfl _l(a)$,
  and hence $m^a_{i,l}>3$ by Lemma~\ref{le:m=2} and by (R4),
  that is $\cm ^a_{i l}\cm ^a_{l i}\ge 2$ by Lemma~\ref{le:Cartan-1}. Similarly
  one gets $\cm ^a_{j l}\cm ^a_{l j}\ge 2$.
  Thus the assumptions in part (2) imply that $t\ge 4$.
\end{proof}

Now we present another technique for the analysis of the finiteness of
$\Wg (\rsC )$. We will use this method in Sect.~\ref{sec:3obj}.

\begin{propo}\label{pr:Gmn}
  Let $m,n\in \ndN $, and define the families $G_{m,n}$ and $H_{m,n}$
  of groups by generators and relations as follows.
  \begin{align*}
    G_{m,n}=&\langle s,t\rangle /(s^2,t^m,(st^{-1}st)^n),\\
    H_{m,n}=&\langle s_1,\ldots ,s_m,T\rangle /
    (s_i^2,T^m, (s_is_{i+1})^n, T^{-1}s_iTs_{i+1}\,|\,1\le i\le m),
  \end{align*}
  where the convention $s_{m+1}=s_1$ is used in the definition of $H_{m,n}$.
  Then there is a group isomorphism $\varphi :G_{m,n}\to H_{m,n}$ with
  $\varphi (s)=s_m$, $\varphi (t)=T$. Further, $G_{m,n}$ is finite if and only
  if $m=1$ or $m=2$ or $n=1$ or $(m,n)=(3,2)$.
\end{propo}

\begin{proof}
  Since $T^{-1}s_mT=s_1$ and $(s_ms_1)^n=1$ in $H_{m,n}$,
  there is a unique group homomorphism $\varphi :G_{m,n}\to H_{m,n}$ with
  $\varphi (s)=s_m$, $\varphi (t)=T$. Further, there is a group homomorphism
  $\psi :H_{m,n}\to G_{m,n}$ with $\psi (s_i)=t^{-i}st^i$, $\psi (T)=t$, and
  the identities $\varphi \psi =\id $ and $\psi \varphi =\id $ hold. Thus
  $\varphi $ is an isomorphism.

  If $m=1$ then $H_{m,n}\simeq \ndZ /2\ndZ $. Assume now that $m\ge 2$.
  Let $N:=\langle s_i\,|\,1\le i\le m\rangle \subset H_{m,n}$. Since
  $T^{-1} s_i T =s_{i+1}$ for all $i\in \{1,2,\ldots ,m\}$, $N$ is a
  normal subgroup of $H_{m,n}$, and $H_{m,n}$ is the semidirect product of $N$
  and the finite abelian group $\langle T\rangle /(T^m)\simeq \ndZ /m\ndZ $.
  Thus $H_{m,n}$ is finite if and only if $N$ is finite. But
  \begin{align*}
    N=\langle s_1,\ldots ,s_m\rangle /(s_i^2,(s_is_{i+1})^n\,|\,1\le i\le m)
  \end{align*}
  is a Coxeter group. It is easy to see that $N$ is finite if and only if
  $n=1$ (and $N\simeq \ndZ /2\ndZ $) or $m=2$ (and $N\simeq \Dih _n$, the
  dihedral group of order $2n$)
  or $(m,n)=(3,2)$ (and $N\simeq (\ndZ /2\ndZ )^3$).
\end{proof}

\begin{propo}\label{pr:Gmnp}
   Let $m,n,p\in \ndN $, and define the families $G_{m,n,p}$ and
   $H_{m,n,p}$
  of groups by generators and relations as follows.
  \begin{align*}
    G_{m,n,p}=&\langle s,u,t\rangle /(s^2,u^2,t^m,(st^{-1}ut)^n,(su)^p),\\
    H_{m,n,p}=&\langle s_1,\ldots ,s_m,u_1,\ldots ,u_m,T\rangle /
    (s_i^2,u_i^2,T^m, (s_iu_{i+1})^n, (s_iu_i)^p,\\
    &\qquad \qquad \qquad T^{-1}s_iTs_{i+1}, T^{-1}u_iTu_{i+1}\,|\,1\le i\le m),
  \end{align*}
  where the convention $s_{m+1}=s_1$, $u_{m+1}=u_1$ is used in the definition
  of $H_{m,n,p}$.
  Then there is a group isomorphism $\varphi :G_{m,n,p}\to H_{m,n,p}$ with
  $\varphi (s)=s_m$, $\varphi (u)=u_m$, $\varphi (t)=T$. Further, $G_{m,n,p}$
  is finite if and only if $m=1$ or $(n,p)=(1,1)$ or $(m,n)=(2,1)$ or
  $(m,p)=(2,1)$ or $(m,n,p)=(3,1,2)$ or $(m,n,p)=(3,2,1)$.
\end{propo}

\begin{proof}
  Entirely similar to the proof of Prop.~\ref{pr:Gmn}.
\end{proof}

\section{Root systems of type $\cC $ with two objects}

Let $I=\{1,2,\ldots ,\theta \}$ for some $\theta \in \ndN $, and $A=\{x,y\}$
with $x\not=y$.
Let $\cC =\cC (I,A,(\rfl _i)_{i\in I},(\Cm ^a)_{a\in A})$ be a connected
Cartan scheme and $\rsC =\rsC (\cC ,(R^a)_{a\in A})$ a root system
of type $\cC $.
Without loss of generality suppose that
\begin{equation}
\begin{aligned}
    \rfl _i(x)=y,\quad \rfl _i(y) =x &\qquad \text{if $1\le i\le \ka $,}\\
    \rfl _i(x)=x,\quad \rfl _i(y) =y &\qquad \text{if $\ka +1\le i\le \theta $}
\end{aligned}
  \label{eq:F2Nact}
\end{equation}
for some $\ka \in I$.
In this case (C2) implies that
$\cm ^x_{i j}=\cm ^y_{i j}$ whenever $1\le i\le \ka $ and $j\in I$.

If $\theta =1$ then $\rfl _1(x)=y$,
$R^x=\{\al _1,-\al _1\}$,
$R^y=\{\al _1,-\al _1\}$, and $\Cm ^x=\Cm ^y=(2)$,
see also \cite[Ex.\,1]{a-HeckYam08}.

Consider now the case $\theta =2$ and $\ka =1$.
Then $\Wg (\rsC )$ is isomorphic to the Coxeter groupoid
\[\langle s_1^x,s_1^y,s_2^x,s_2^y\rangle
/(s_1^y s_1^x, s_1^x s_1^y, (s_2^x)^2, (s_2^y)^2, (s_1^y s_2^y s_1^x s_2^x)^m)
\]
for some $m\in \ndN $, see Thm.~\ref{th:Coxgr}.
Identify $\s_i^a$, where $i\in \{1,2\}$ and $a\in \{x,y\}$,
with its matrix with respect to the standard basis
$\{\al _1,\al _2\}$ of $\ndZ ^2$. Further, note that
$\cm ^x_{1 2}=\cm ^y_{1 2}$ by (C2). One gets
\begin{align*}
  \s _1\s _2\s _1\s _2^x
  &=
  \begin{pmatrix}
    -1 & -\cm ^y_{1 2} \\ 0 & 1
  \end{pmatrix}
  \begin{pmatrix}
    1 & 0 \\ -\cm ^y_{2 1} & -1
  \end{pmatrix}
  \begin{pmatrix}
    -1 & -\cm ^x_{1 2} \\ 0 & 1
  \end{pmatrix}
  \begin{pmatrix}
    1 & 0 \\ -\cm ^x_{2 1} & -1
  \end{pmatrix}\\
  &=
  \begin{pmatrix}
    \cm ^x_{1 2}\cm ^y_{2 1}-1 & \cm ^x_{1 2} \\ -\cm ^y_{2 1} & -1
  \end{pmatrix}
  \begin{pmatrix}
    \cm ^x_{1 2}\cm ^x_{2 1}-1 & \cm ^x_{1 2} \\ -\cm ^x_{2 1} & -1
  \end{pmatrix}\\
  &=
  \begin{pmatrix}
    \cm ^x_{1 2}\cm ^x_{2 1}\cm ^x_{1 2}\cm ^y_{2 1}-\cm ^x_{1 2}\cm ^y_{2 1}
    -2\cm ^x_{1 2}\cm ^x_{2 1}+1 &
    \cm ^x_{1 2}(\cm ^x_{1 2}\cm ^y_{2 1}-2) \\
    -\cm ^x_{1 2}\cm ^x_{2 1}\cm ^y_{2 1}+\cm ^y_{2 1}+\cm ^x_{2 1} &
    1-\cm ^x_{1 2}\cm ^y_{2 1}
  \end{pmatrix}.
\end{align*}
The groupoid $\Wg (\rsC )$ is finite if and only if
$\s :=\s _1\s _2\s _1\s _2^x$ has finite order.
Since $\s \in \Aut (\ndZ ^2)$, the latter is equivalent to the condition
\begin{align}
  \s =\pm \id \qquad \text{or} \qquad \tr \s \in \{-1,0,1\}.
  \label{eq:sfin}
\end{align}
By (M2) and (C2) one has
$\s =\id $ if and only if
$\cm ^x_{1 2}=\cm ^x_{2 1}=\cm ^y_{1 2}=\cm ^y_{2 1}=0$.
Further, $\s =-\id $ if and only if
$\cm ^x_{1 2}\cm ^y_{2 1}=2$, $\cm ^x_{2 1}=\cm ^y_{2 1}$, that is
\begin{align*}
  \Cm ^x=\Cm ^y=
  \begin{pmatrix}
    2 & -1 \\ -2 & 2
  \end{pmatrix}
  \qquad \text{or} \qquad
  \Cm ^x=\Cm ^y=
  \begin{pmatrix}
    2 & -2 \\ -1 & 2
  \end{pmatrix}.
\end{align*}
It remains to consider the second relation in Eq.~\eqref{eq:sfin}. One has
\begin{align*}
  \tr \s =(\cm ^x_{1 2}\cm ^x_{2 1}-2)(\cm ^x_{1 2}\cm ^y_{2 1}-2)-2,
\end{align*}
and hence relation $-1\le \tr \s \le 1$ is equivalent to
\begin{align*}
  1\le (\cm ^x_{1 2}\cm ^x_{2 1}-2)(\cm ^x_{1 2}\cm ^y_{2 1}-2)\le 3.
\end{align*}
By (M1) and (M2) one gets
$\cm ^x_{1 2}=\cm ^x_{2 1}=\cm ^y_{1 2}=\cm ^y_{2 1}=-1$ or
$\cm ^x_{1 2}\cm ^x_{2 1}\ge 3$,
$\cm ^x_{1 2}\cm ^y_{2 1}\ge 3$.
However, the first case contradicts (R4) and Lemma~\ref{le:Cartan-1}.
Now the following statement can be obtained easily.

\begin{propo}\label{pr:r=2,k=1}
  Assume that $\theta =2$ and $\ka =1$. Then $\rsC $ is finite if and only if
  the Cartan matrices $\Cm ^x$, $\Cm ^y$ satisfy, up to permutation of $x$ and
  $y$, one of the following conditions (1),(2).
  \begin{enumerate}
    \item $\Cm ^x=\Cm ^y$ is of finite type $A_1\times A_1$, $B_2$, or $G_2$,
      that is, $\cm ^x_{1 2}=\cm ^x_{2 1}=0$
      or $\cm ^x_{1 2}\cm ^x_{2 1}\in \{2,3\}$.
    \item $\cm ^x_{1 2}=\cm ^y_{1 2}=-1$, $\cm ^x_{2 1}=-3$,
      and $\cm ^y_{2 1}\in \{-4,-5\}$.
  \end{enumerate}
  In case (1) $R^x=R^y$ is the usual set of roots
  corresponding to the generalized Cartan matrix $\Cm ^x=\Cm ^y$.
  In case (2) one has
  \begin{align*}
    R^x_+=&\{1,2,12,12^2,12^3,1^2 2^3,1^3 2^4,1^3 2^5\},\\
    R^y_+=&\{1,2,12,12^2,12^3,1 2^4,1^2 2^3,1^2 2^5\}
  \end{align*}
  if $\cm ^y_{2 1}=-4$, and
  \begin{align*}
    R^x_+=&\{1,2,12,12^2,1 2^3,1^2 2^3,1^3 2^4, 1^3 2^5,
    1^4 2^5, 1^4 2^7, 1^5 2^7, 1^5 2^8 \},\\
    R^y_+=&\{1,2,12,12^2,12^3,12^4,12^5,1^2 2^3,1^2 2^5,
    1^2 2^7, 1^3 2^7, 1^3 2^8\}
  \end{align*}
  if $\cm ^y_{2 1}=-5$.
\end{propo}

In the last part of the proposition the abbreviation
$1^m 2^n=m\al _1+n\al _2$ was used, where exponents $1$ and factors
$i^0$ for $i\in \{1,2\}$ are omitted.
The determination of $R_+^a$ is straightforward, see
Prop.~\ref{pr:allposroots} or, more directly, \cite[Lemma\,6]{a-HeckYam08}.

The case $\theta=2$, $\ka =2$ is even easier than the case $\ka =1$.
Indeed, by (C2) one has
$\Cm ^x=\Cm ^y$, and hence Prop.~\ref{pr:allposroots} implies that
$R^x=R^y$ is the root system of rank $2$
corresponding to the Cartan matrix $\Cm ^x$. Thus the following holds.

\begin{propo}\label{pr:r=2,k=2}
  Assume that $\theta =2$ and $\ka =2$. Then $\rsC $ is finite if and only if
  $\rsC $ is standard and $\Cm ^x=\Cm ^y$ is of finite type,
  that is,
  $\cm ^x_{1 2}=\cm ^x_{2 1}=0$ or $\cm ^x_{1 2}\cm ^x_{2 1}\in \{1,2,3\}$.
\end{propo}

Now assume that $\theta \ge 3$. As before, suppose
that $1\le \ka \le \theta $ such that Eq.~\eqref{eq:F2Nact} holds.
First we develop some general properties.

\begin{lemma}
  Suppose that $\rsC $ is finite and that $1\le \ka \le \theta -2$.
  Then $\cm ^x_{l i}\cm ^x_{l j}=0$
  for all $l\le \ka $ and $i,j>\ka $ with $i\not=j$.
  \label{le:1connect}
\end{lemma}

\begin{proof}
  Let $I'=\{l,i,j\}$, and consider the restriction $\rsC |_{I'}$ of $\rsC $,
  see Sect.~\ref{sec:dec}. Since $\rfl _l(x)=y$, $\rsC |_{I'}$ is
  connected. Thus it suffices to consider the case $I=I'$.
  Let $t_i:=\s _l^y \s _i^y \s _l^x$, $t_j:=\s _l^y \s _j^y \s _l^x$.
  By Prop.~\ref{pr:Stabrel},
  the subgroup $\Stab (x)$ of $\Hom (\Wg (\rsC ))$
  is generated by $s_i:=\s _i^x$,
  $s_j:=\s _j^x$, $t_i$, and $t_j$. Moreover, $\Stab (x)$
  can be presented by the relations
  \begin{align*}
    s_i^2,s_j^2,t_i^2,t_j^2, (s_is_j)^{m_{i,j}^x}, (t_i t_j)^{m_{i,j}^y},
    (s_i t_i)^{m_{l,i}^x/2}, (s_j t_j)^{m_{l,j}^x/2}.
  \end{align*}
  Assume now that $\cm ^x_{l i}\cm ^x_{l j}\not=0$.
  Then $m_{l,i}^x>2$ by Lemma~\ref{le:m=2}, that is $m_{l,i}^x/2\ge 2 $, and similarly
  $m_{l,j}^x/2\ge 2$.
  Thus the subgroup of $\Stab (x)$
  generated by $s_i$ and $t_j$ is infinite and hence $\Stab (x)$ is an
  infinite Coxeter group. This is a contradiction to the finiteness of $\rsC $,
  and hence the claim of the lemma is proven.
\end{proof}

\begin{theor}\label{th:o=2}
  Let $\theta \in \ndN $, $I=\{1,2,\ldots ,\theta \}$, $A=\{x,y\}$,
  and $\ka \in I$ as in
  Eq.~\eqref{eq:F2Nact}. If $\rsC $ is finite and irreducible, then
  up to permutation of $\{x,y\}$ and $I$,
  one of the following sets of conditions hold.
\begin{enumerate}
  \item $\theta \in \ndN $, $\ka \in \{1,2,\ldots ,\theta \}$,
    and $\Cm ^x=\Cm ^y$ is an indecomposable
    Cartan matrix of finite type such that $\cm ^x_{i j}\cm ^x_{j i}\not=1$
    for all $i,j$ with $i\le \ka $, $j>\ka $.
  \item $\theta =2$, $\ka =1$, $\cm ^x_{1 2}=\cm ^y_{1 2}=-1$, $\cm ^x_{2 1}=-3$,
    $\cm ^y_{2 1}\in \{-4,-5\}$.
  \item $\theta =3$, $\ka =1$,
    \begin{align*}
      \Cm ^x=
      \begin{pmatrix}
        2 & -1 & 0 \\ -2 & 2 & -1 \\ 0 & -1 & 2
      \end{pmatrix},\quad
      \Cm ^y=
      \begin{pmatrix}
        2 & -1 & 0 \\ -2 & 2 & -2 \\ 0 & -1 & 2
      \end{pmatrix}.
    \end{align*}
  \item $\theta =3$, $\ka =1$,
     \begin{align*}
      \Cm ^x=
      \begin{pmatrix}
        2 & -2 & 0 \\ -1 & 2 & -1 \\ 0 & -1 & 2
      \end{pmatrix},\quad
      \Cm ^y=
      \begin{pmatrix}
        2 & -2 & 0 \\ -1 & 2 & -2 \\ 0 & -1 & 2
      \end{pmatrix}.
    \end{align*}
\end{enumerate}
Conversely, if $\Cm ^x,\Cm ^y$ satisfy one of the conditions (1)--(4), then
$\rsC $ is finite and irreducible.
\end{theor}

\begin{proof}
  If $\theta \le 2$ then the claim of the theorem holds by
  Props.~\ref{pr:r=2,k=1}, \ref{pr:r=2,k=2}.
  
  Assume that $\theta \ge 3$. Consider first the case $\ka =1$,
  that is, $\rfl _1(x)=y$, $\rfl _1(y)=x$, and
  $\rfl _i(x)=x$, $\rfl _i(y)=y$ if $2\le i\le \theta $.
  Using a permutation of $\{2,3,\ldots ,\theta \}$,
  by Lemma~\ref{le:1connect} one can assume that $\cm ^x_{1 i}=0$
  for all $i\ge 3$.
  Then
  \begin{align}
    \cm ^x_{i l}=\cm ^y_{i l} \qquad \text{for all $i\ge 3$ and $l\in I$}
    \label{eq:c=c}
  \end{align}
  by Lemma~\ref{le:samerowfor0}, and hence
  $\cm ^x_{i 1}=\cm ^y_{1 i}=\cm ^y_{i 1}=0$ for all $i\ge 3$ by (M2).
  We obtain that
  $\s _1^y\s _i^y\s _1^x=\s _i^x$ for all $i\ge 3$ by Thm.~\ref{th:Coxgr}.
  Therefore Prop.~\ref{pr:Stabrel} tells that
  $\Hom (x)$ is isomorphic to the group 
  \begin{equation}
    \begin{split}
    \langle s_2,s_3,\ldots ,s_\theta ,t\rangle /&
    ( s_i^2,t^2, (s_j s_k)^{m_{j,k}^x},(t s_l)^{m_{2,l}^y},(s_2t)^{m_{1,2}^x/2}
    \,|\\
    &2\le i\le \theta , 2\le j<k\le \theta , 2<l\le \theta ),
  \end{split}
    \label{eq:k=1stab}
  \end{equation}
  where $t$ corresponds to the element $\s _1^y\s _2^y\s _1^x$.
  Since $\rsC $ is irreducible and $m_{1,i}^x=2$ for $i\ge 3$,
  we obtain that $m_{1,2}^x\ge 3$ by Lemma~\ref{le:m=2} and
  Prop.~\ref{pr:decomp}.
  Thus $\Hom (x)$ is a Coxeter group of rank $\theta $.
  Further, again since $\rsC $ is irreducible and
  $\cm ^x_{1 i}=0$ for all
  $i\ge 3$, we have $\cm ^x_{2 r}\not=0$ for some $r\ge 3$. Without loss of
  generality assume that $\cm ^x_{2 3}\not=0$.
  Then $\cm ^y_{3 2}=\cm ^x_{3 2}\not=0$ by Eq.~\eqref{eq:c=c} and (M2),
  and hence $\cm ^y_{2 3}\not=0$ and $m_{2,3}^x,m_{2,3}^y\ge 3$.
  Since $\Hom (x)$ is a finite Coxeter group,
  this implies that $s_2t=ts_2$, that is,
  $m_{1,2}^x=4$. Then Prop.~\ref{pr:r=2,k=1} gives that
  $\cm ^x_{i j}=\cm ^y_{i j}$ for $i,j\in \{1,2\}$,
  and $\cm ^x_{1 2}\cm ^x_{2 1}=2$.
  
  First assume that $\theta \ge 4$ and
  $\cm ^x_{2 i}\not=0$ for some $i\ge 4$. As above, we
  conclude that $m_{2,i}^x,m_{2,i}^y\ge 3$, which is a contradiction to the
  finiteness of $\Hom (x)$ and the relations
  $m_{2,3}^x,m_{2,3}^y\ge 3$. Hence, if $\theta \ge 4$, then
  $\cm ^x_{2 i}=\cm ^y_{2 i}=0$ for all $i\ge 4$.
  Since $\rsC $ is irreducible, there exists
  $r\ge 4$ with $\cm ^x_{3 r}\not=0$. Then $m_{2,3}^x=m_{2,3}^y=3$
  by the finiteness of the group in Eq.~\eqref{eq:k=1stab}, and hence
  $\cm ^x_{2 i}=\cm ^y_{2 i}$ for all $i$. Thus, if $\theta \ge 4$, then
  $\Cm ^x=\Cm ^y$ as in (1).
  The restriction $\cm^x_{1 j}\cm ^x_{j 1}\not=1$ comes
  from Lemma~\ref{le:Cartan-1} and
  the relation $\rfl _1\rfl _j\rfl _1(x)=x\not=y=\rfl _j\rfl _1\rfl _j(x)$
  for $j>1$.

  If $\theta =3$ and $\ka =1$, then Eq.~\eqref{eq:k=1stab} tells that
  \begin{equation}
    \Hom (x)\simeq \langle s_2,s_3,t\rangle /
    ( s_2^2,s_3^2,t^2, (s_2 s_3)^{m_{2,3}^x},(t s_3)^{m_{2,3}^y},s_2ts_2t).
    \label{eq:r=3,k=1stab}
  \end{equation}
  Further, $m_{2,3}^x,m_{2,3}^y\ge 3$ by the above arguments. Thus,
  $(m^x_{2,3},m^y_{2,3})\in \{(3,3),(3,4),(4,3)\}$.
  Since $\rfl _2(a)=\rfl _3(a)=a$ for $a\in \{x,y\}$, Thm.~\ref{th:standard}
  yields that
  $(\cm ^x_{2 3}\cm ^x_{3 2},\cm ^y_{2 3}\cm ^y_{3 2})\in
  \{(1,1),(1,2),(2,1)\}$.
  Recall that $\cm ^x_{i j}=\cm ^y_{i j}$ for all $(i,j)\not=(2,3)$.
  Thus, if $\cm ^a_{2 3}\cm ^a_{3 2}=1$ for $a\in \{x,y\}$, then $\rsC $
  satisfies the conditions in (1). Otherwise, since $\cm ^x_{3 2}=\cm ^y_{3
  2}$, $\rsC $ satisfies up to permutation of $x$ and $y$ the conditions of (3)
  or (4).

  Consider now the case $\theta \ge 3$, $\ka \ge 2$. From (C2)
  we obtain that $\cm ^x_{i n}=\cm ^y_{i n}$ for all $i\le \ka $ and $n\in I$.
  Further, Lemma~\ref{le:Cartan3}(2) for $l>\ka $, $i,j\le \ka $, $a:=x$, and $b:=y$
  implies that $\cm ^x_{l n}=\cm ^y_{l n}$ for all $l>\ka $ and $n\in I$.
  Thus, $\rsC $ satisfies the conditions in (1), where the
  restriction $\cm^x_{i j}\cm ^x_{j i}\not=1$ comes again
  from Lemma~\ref{le:Cartan-1} and
  the relation $\rfl _i\rfl _j\rfl _i(x)=x\not=y=\rfl _j\rfl _i\rfl _j(x)$
  for $i\le \ka $, $j>\ka $.
\end{proof}

\begin{remar}\label{re:WgNichols2}
  The appearance of non-standard root systems in Thm.~\ref{th:o=2}
  is not surprizing. With an appropriate definition of the Weyl groupoid of a
  Nichols algebra of diagonal type, the examples in Thm.~\ref{th:o=2}(2) can
  be identified with the root systems of the Nichols algebras corresponding
  to Row~14 and Row~17 of \cite[Table\,1]{a-Heck05b}, respectively. Similarly,
  the examples in Thm.~\ref{th:o=2}(3),(4) can
  be identified with the root systems of the Nichols algebras corresponding
  to Row~13 and Row~18 of \cite[Table\,2]{a-Heck05b}, respectively.
\end{remar}

\section{Root systems of type $\cC $ with three objects}
\label{sec:3obj}

Let $\theta \in \ndN $, $I=\{1,2,\ldots ,\theta \}$, and $A$ a set of
cardinality $3$.
Let $\cC =\cC (I,A,(\rfl _i)_{i\in I},(\Cm ^a)_{a\in A})$ be a connected
Cartan scheme and $\rsC =\rsC (\cC ,(R^a)_{a\in A})$ a root system
of type $\cC $.
In this case we necessarily have $\theta \ge 2$.

Let first $\theta =2$.
Then, up to enumeration of the objects and up to permutation of $I$,
we may fix the three elements $x,y$, and $z$ of $A$, such that
the object change diagram of $\rsC $ is
\setlength{\unitlength}{1mm}
\begin{align}
  \begin{picture}(26,10)
    \put(3,3){\circle*{2}}
    \put(4.5,3){\line(1,0){7}}
    \put(13,3){\circle*{2}}
    \put(14.5,3){\line(1,0){7}}
    \put(23,3){\circle*{2}}
    \put(2,6){$x$}
    \put(12,6){$y$}
    \put(22,6){$z$}
    \put(7,5){$_1$}
    \put(17,5){$_2$}
  \end{picture}
  \label{eq:obchanger2}
\end{align}

By (C2) the generalized Cartan matrices are
\begin{align}
  \Cm ^x=
  \begin{pmatrix}
    2 & -a \\ -c & 2
  \end{pmatrix},\quad
  \Cm ^y=
  \begin{pmatrix}
    2 & -a \\ -d & 2
  \end{pmatrix},\quad
  \Cm ^z=
  \begin{pmatrix}
    2 & -b \\ -d & 2
  \end{pmatrix},
  \label{eq:r2o3cm}
\end{align}
where $a,b,c,d\in \ndN _0$. 
Further --- by replacing $(1,2)$ by $(2,1)$ and $(x,z)$ by $(z,x)$,
if necessary ---
one may assume that $a\le d$.

\begin{theor}\label{th:r2o3}
  Let $\rsC $ be a connected root system of type $\cC $ of rank $2$ with $3$
  objects.
  Assume that $A=\{x,y,z\}$ and $I=\{1,2\}$ such that the object change
  diagram of $\rsC $ is as in Eq.~\eqref{eq:obchanger2}. Further, let
  $a,b,c,d\in \ndN _0$ such that $a\le d$ and that
  the Cartan matrices $\Cm ^x,\Cm ^y$ and $\Cm ^z$ are as in
  Eq.~\eqref{eq:r2o3cm}. If $\rsC $ is finite, then
  $(a,b,c,d)$ and $R^x$ satisfy one of the
  following equations.
\begin{enumerate}
  \item $(a,b,c,d)=(1,1,1,1)$, $|R^x_+|=3$.
  \item $(a,b,c,d)=(1,1,3,3)$, $|R^x_+|=6$.
  \item $(a,b,c,d)=(1,2,4,2)$, $|R^x_+|=6$.
  \item $(a,b,c,d)=(1,3,6,2)$, $|R^x_+|=12$.
  \item $(a,b,c,d)=(1,4,5,2)$, $|R^x_+|=12$.
  \item $(a,b,c,d)=(1,3,7,2)$, $|R^x_+|=18$.
  \item $(a,b,c,d)=(1,5,5,2)$, $|R^x_+|=18$.
\end{enumerate}
Conversely, if $(a,b,c,d)$ is one of the above $7$ quadruples, then $\rsC $ is
finite.
\end{theor}

\begin{proof}
If one of the relations $a=0$,
$b=0$, $c=0$, and $d=0$ holds, then also the other three because of
(M2). However,
then $x=\rfl _2\rfl _1(y)=\rfl _1\rfl _2(y)=z$ by
(R4) and Lemma~\ref{le:m=2}, which is a contradiction. Thus $a,b,c,d>0$.

Clearly, the finiteness of $\rsC $ implies that the order of the linear map
$\tilde{\s }:=\s _2^x\s _1^y\s _2^z\s _1^z\s _2^y\s _1^x\in \Aut (\ndZ ^2)$
is finite, that is, the matrix $t=(t_{ij})_{i,j=1,2}$, where
\begin{align*}
  t_{11}=& -abd^2+2ad+bd-1,\\
  t_{12}=& a^2bd^2-2a^2d-2abd+2a+b,\\
  t_{21}=& -abcd^2+2acd+bcd+bd^2-c-2d,\\
  t_{22}=& a^2bcd^2-2a^2cd-2abcd-abd^2+2ac+2ad+bc+bd-1,
\end{align*}
has finite order. This means that either $t=\id $ or $t=-\id $ or
$t_{11}+t_{22}\in \{-1,0,1\}$.
Observe that
\begin{align}
  t_{11}+t_{22}-2=&(ad-1)(abcd-2ac-bc-2bd+4),\label{eq:t1}\\
  t_{11}+t_{22}+2=&(abd-2a-b)(acd-c-2d),\label{eq:t2}\\
  t_{22}=&(1-ac)t_{11}+c(-abd+a+b).\label{eq:t3}
\end{align}
Thus, if $t=\id $, then $1-ac-abcd+ac+bc=1$ by Eq.~\eqref{eq:t3},
that is, $ad=1$.
Therefore $a=d=1$, and relation $t_{21}=0$ gives that $b+c=2$.
Since $bc\not=0$, we obtain that $b=c=1$. This gives solution (1).

Suppose now that $t=-\id $. Then $-1+2ac-abcd+bc=-1$ by Eq.~\eqref{eq:t3},
that is, $2a+b=abd$. Inserting this into Eq.~\eqref{eq:t1} we get that
$-4=(ad-1)(4-2bd)$. Since $d\ge a$, we obtain the solutions $(a,b,d)\in
\{(1,2,2),(1,1,3)\}$. Using equation $t_{21}=0$ we get solutions (2) and (3).

It remains to determine all solutions with $t_{11}+t_{22}\in \{-1,0,1\}$.
Eq.~\eqref{eq:t1} yields that $ad\in \{2,3,4\}$. If $ad=4$, then
either $(a,d)=(1,4)$ or $(a,d)=(2,2)$. Further, by
Eq.~\eqref{eq:t2} we get that
\[(3b-2a)(3c-2d)\in \{1,2,3\}.\]
Both for $(a,d)=(1,4)$ and for
$(a,d)=(2,2)$ there is a unique solution of this relation with $b,c\in \ndN $,
namely, $(a,b,c,d)=(1,1,3,4)$ and $(a,b,c,d)=(2,1,1,2)$.
However, in the first case one has
\[ (\s _1\s _2)^4\s_1^x(\al _1)=\al _1-\al _2,\]
and in the second case one gets $\s _2\s _1\s _2\s _1^x
(\al _2)=\al _1-\al _2$. Therefore (R1)
gives that there is no root system of type $\cC $ with $ad=4$.

Assume now that $ad=3$, that is, $a=1$ and $d=3$.
Then Eq.~\eqref{eq:t1} implies that $-3\le 2(2bc-2c-6b+4)\le -1$, which has no
solution with $a,b,c,d\in \ndN $.

Finally, let $a=1$ and $d=2$. Then, by Eq.~\eqref{eq:t2}, $t_{11}+t_{22}\in
\{-1,0,1\}$ if and only if
\begin{align}\label{eq:bcrel}
  1\le (b-2)(c-4)\le 3.
\end{align}
If $b=1$ then $c\in \{1,2,3\}$. If $(b,c)=(1,1)$ then
$\s _2^x\s _1^y\s _2^z(\al _1)=\al _1-\al _2$, a contradiction to (R1).
If $(b,c)=(1,2)$ then $\Cm ^x=\Cm ^y=\Cm ^z$, $|R_+^x|=4$,
but $(\rfl _1\rfl _2)^4(x)=y\not=x$, a contradiction to (R4).
If $(b,c)=(1,3)$ then $\s _1^z\s _2^y\s _1^x\s _2^x(\al _1)
=\al _2-\al _1$, a contradiction to (R1).

If $b=3$ in Rel.~\eqref{eq:bcrel}, then $c\in \{5,6,7\}$.
If $c=5$ then direct computation shows that
$(\s _1\s _2)^4 1_y(\al _1)=\al _2-\al _1$,
a contradiction to (R1).
If $c=6$, then we get solution~(4), and if $c=7$, then we get solution~(6).

The remaining two solutions of Rel.~\eqref{eq:bcrel} are $(b,c)=(4,5)$ and
$(b,c)=(5,5)$. These correspond to solution (5) and (7), respectively.

The sets $R_+^x$ can be calculated from Prop.~\ref{pr:allposroots}.
\end{proof}

\begin{remar}\label{re:WgNichols3}
  It is interesting to note that in contrast to the case with two objects, 
  see Rem.~\ref{re:WgNichols2}, not all non-standard root systems
  can be obtained from Nichols algebras of diagonal type.
  The example in Thm.~\ref{th:r2o3}(3) can
  be identified with the root system of the Nichols algebras corresponding
  to Row~10 \cite[Table\,1]{a-Heck05b}.
  However, the only rank two Nichols algebras of diagonal type with
  12 positive roots are those in Row~17 of \cite[Table\,1]{a-Heck05b},
  and there are no such Nichols algebras with more than 12 positive roots.
  This can be read off from the trees in the appendix of \cite{a-Heck07a}.
  The Nichols algebras corresponding to
  Row~17 of \cite[Table\,1]{a-Heck05b} have been discussed already in
  Rem.~\ref{re:WgNichols2}: in all of the Cartan matrices one has at least
  one entry $-1$. Thus the examples in Thm.~\ref{th:r2o3}(4)--(7) can not be
  obtained as the root system of a Nichols algebra of diagonal type.

  It is not clear if there are more general Nichols algebras with such root
  systems.
\end{remar}

Now we assume that $\theta =3$.

\begin{theor}
  \label{th:o=3,r=3}
  Let $\cC $ be a connected Cartan scheme with $I=\{1,2,3\}$ and with $3$
  objects, and let
  $\rsC $ be a finite irreducible root system of type $\cC $.
  Then $\rsC $ is standard, and the Cartan matrices
  are indecomposable and of type $A_3$, $B_3$, or
  $C_3$. If $\cm ^a_{12}=\cm ^a_{21}=-1$ and
  $\cm ^a_{13}=\cm ^a_{31}=0$ for all objects $a$,
  then the object change diagram of $\rsC $ is
  \setlength{\unitlength}{1mm}
  \begin{align*}
    &
    \begin{picture}(26,6)(0,2)
      \put(3,3){\circle*{2}}
      \put(4.5,3.5){\line(1,0){7}}
      \put(4.5,2.5){\line(1,0){7}}
      \put(13,3){\circle*{2}}
      \put(14.5,3){\line(1,0){7}}
      \put(23,3){\circle*{2}}
      \put(8,5){\makebox[0pt]{$_{1,3}$}}
      \put(18,5){\makebox[0pt]{$_2$}}
    \end{picture}
    &&\text{for type $A_3$ Cartan matrices, and}\\&
    \begin{picture}(26,6)(0,2)
      \put(3,3){\circle*{2}}
      \put(4.5,3){\line(1,0){7}}
      \put(13,3){\circle*{2}}
      \put(14.5,3){\line(1,0){7}}
      \put(23,3){\circle*{2}}
      \put(8,5){\makebox[0pt]{$_1$}}
      \put(18,5){\makebox[0pt]{$_2$}}
    \end{picture}
    &&\text{for type $B_3$ and $C_3$ Cartan matrices.}
  \end{align*}
\end{theor}

\begin{proof}
  Let $x,y,z$ denote the three objects of $\rsC $. Since $\cC $ is connected,
  we may assume (using permutations of $I$ and $A$) that the
  restriction of $\rsC $ to $I=\{1,2\}$ is as for $\theta =2$ --- without
  supposing that $a\le d$ in Eq.~\eqref{eq:r2o3cm}.
  Then we have to consider three cases.
  
  \textit{Case 1:} Assume that $\rfl _3:A\to A$ is the identity.
  By
  Thm.~\ref{th:Coxgr} and Prop.~\ref{pr:Stabrel} the group $\Hom (x)$ is
  isomorphic to
  \begin{equation}
  \begin{aligned}
    \Hom (x)\simeq \langle &s_2,s_3,t_3,u_1,u_3\rangle
    /(s_2^2,s_3^2,t_3^2,u_1^2,u_3^2,
    (u_1s_2)^{m_{1,2}^x/3},\\
    &(t_3s_3)^{m_{1,3}^x/2},
    (s_2s_3)^{m_{2,3}^x}, (u_3 t_3)^{m_{2,3}^y/2}, (u_1u_3)^{m_{1,3}^z}),
  \end{aligned}
    \label{eq:Homx}
  \end{equation}
  where the isomorphism is given by $\s _2^x\mapsto s_2$, $\s _3^x\mapsto
  s_3$, $\s _1^y\s _3^y\s _1^x\mapsto t_3$,
  $\s _1^y\s _2^z\s _1^z\s _2^y\s _1^x\mapsto u_1$, and
  $\s _1^y\s _2^z\s _3^z\s _2^y\s _1^x\mapsto u_3$.
  Thus, $\Hom (x)$ is a Coxeter group.
  
  Suppose first that $m_{1,3}^x=m_{2,3}^y=2$. Then
  $\cm ^x_{1 3}=\cm ^y_{2 3}=0$, and hence $\cm ^y_{1 3}=0$ by
  (C2) and relation $\rfl _1(x)=y$.
  Thus, $\rsC $ is not irreducible by
  Prop.~\ref{pr:decomp}, which is a contradiction. Note that
  $m_{2,3}^z=m_{2,3}^y$, since $\rfl _2(y)=z$. Now, using a symmetry in
  the presentation of $\rsC $, we can assume that $m_{2,3}^y>2$.

  If $m_{1,2}^x>3$, then the quotient
  \begin{align*}
    \Hom (x)/(t_3s_3)\simeq \langle s_2,s_3,&u_1,u_3\rangle
    /(s_2^2,s_3^2,u_1^2,u_3^2,
    (u_1s_2)^{m_{1,2}^x/3},\\
    &(s_2s_3)^{m_{2,3}^x}, (u_3 s_3)^{m_{2,3}^y/2}, (u_1u_3)^{m_{1,3}^z})
  \end{align*}
  is a Coxeter group without relation between $u_1$ and $s_3$, and hence it is
  infinite, which is a contradiction to the finiteness of $\rsC $.
  If $m_{1,2}^x=3$, then
  \begin{align*}
    \Hom (x)\simeq \langle &s_2,s_3,t_3,u_3\rangle
    /(s_2^2,s_3^2,t_3^2,u_3^2,\\
    &(t_3s_3)^{m_{1,3}^x/2},
    (s_2s_3)^{m_{2,3}^x}, (u_3 t_3)^{m_{2,3}^y/2}, (s_2u_3)^{m_{1,3}^z}).
  \end{align*}
  In this case, if $m_{1,3}^x>2$, then there is no Coxeter relation between
  $s_2$ and $t_3$, which is again a contradiction. Thus $m_{1,3}^x=2$ and
  \begin{equation}
  \begin{aligned}
    \Hom (x)\simeq \langle &s_2,s_3,u_3\rangle
    /(s_2^2,s_3^2,u_3^2,\\
    &(s_2s_3)^{m_{2,3}^x}, (u_3 s_3)^{m_{2,3}^y/2}, (s_2u_3)^{m_{1,3}^z}).
  \end{aligned}
    \label{eq:Homx1}
  \end{equation}
  Since $m_{1,2}^x=3$, $m_{1,3}^x=2$ and $m_{2,3}^y>2$, Thm.~\ref{th:r2o3}
  yields that the Cartan matrices
  $\Cm ^x$, $\Cm ^y$, and $\Cm ^z$ are
\begin{align*}
  \Cm ^x=
  \begin{pmatrix}
    2 & -1 & 0 \\
    -1 & 2 & -a \\
    0 & -b & 2
  \end{pmatrix},
  \Cm ^y=
  \begin{pmatrix}
    2 & -1 & 0 \\
    -1 & 2 & -c \\
    0 & -d & 2
  \end{pmatrix},
  \Cm ^z=
  \begin{pmatrix}
    2 & -1 & -e \\
    -1 & 2 & -c \\
    -f & -g & 2
  \end{pmatrix},
\end{align*}
where $a,b,e,f,g\in \ndN _0$ and $c,d\in \ndN $.
Moreover, Lemmata~\ref{le:samerowfor0}, \ref{le:m=3} imply that
\begin{align}
  d=b,\qquad a=c+e.
  \label{eq:acrel}
\end{align}
The isomorphism in Eq.~\eqref{eq:Homx1} tells that $\s _2^x$ and $\s _3^x$
generate a finite Coxeter subgroup of $\Hom (x)$, and hence $ab\in
\{0,1,2,3\}$. Since $a=c+e$ and $c>0$, we get $ab\in \{1,2,3\}$, and hence
$m_{2,3}^x\ge 3$.

If $ab=3$, then $m_{2,3}^x=6$, and therefore the finiteness
of $\Hom (x)$ and Eq.~\eqref{eq:Homx1} imply that $m_{2,3}^y=4$,
$m_{1,3}^z=2$.
The former gives by Prop.~\ref{pr:r=2,k=1} that $cd=2$, and the latter, that
$e=f=0$. By Eq.~\eqref{eq:acrel} we get $2=cd=ab=3$, a contradiction.

If $ab=1$, then $a=b=1$, and hence Eq.~\eqref{eq:acrel} and relation $c\in \ndN
$ imply that $e=0$, $c=a=1$, and $d=b=1$. Thus $m_{2,3}^y=3$ by
Lemma~\ref{le:Cartan-1}, which is a contradiction to
$\rfl _2\rfl _3\rfl _2(y)=y\not=z=\rfl _3\rfl _2\rfl _3(y)$ and (R4).

Finally, assume that $ab=2$. If $e>0$, then $a=2$, $b=1$, $c=1$, $d=1$, and
$e=1$ by Eq.~\eqref{eq:acrel}. Then $m_{2,3}^y=3$ by Lemma~\ref{le:Cartan-1},
but $\rfl _2\rfl _3\rfl _2(y)\not=\rfl _3\rfl _2\rfl _3(y)$,
which contradicts (R4).
Thus $e=0$, and hence $c=a$, $d=b$, and $f=0$ by Eq.~\eqref{eq:acrel} and
(M2).
Since $cd=2$, we get $g=d$ by
Prop.~\ref{pr:r=2,k=1}, and hence $\Cm ^x=\Cm ^y=\Cm ^z$ are Cartan matrices
of type $B_3$ or $C_3$.

\textit{Case 2:} Assume that $\rfl _3(x)=y$, $\rfl _3(y)=x$, and $\rfl _3(z)=z$.
In other words, the object change diagram of $\rsC $ is
\setlength{\unitlength}{1mm}
\begin{align}
  \begin{picture}(26,6)(0,2)
    \put(3,3){\circle*{2}}
    \put(4.5,3.5){\line(1,0){7}}
    \put(4.5,2.5){\line(1,0){7}}
    \put(13,3){\circle*{2}}
    \put(14.5,3){\line(1,0){7}}
    \put(23,3){\circle*{2}}
    \put(3,6){\makebox[0pt]{$x$}}
    \put(13,6){\makebox[0pt]{$y$}}
    \put(23,6){\makebox[0pt]{$z$}}
    \put(8,5){\makebox[0pt]{$_{1,3}$}}
    \put(18,5){\makebox[0pt]{$_2$}}
  \end{picture} .
  \label{eq:obchanger3}
\end{align}
%
We apply Prop.~\ref{pr:Stabrel} to
$\Hom (\rsC )$. Let $X_x=1_x$, $X_y=\s _1^x\in \Hom (x,y)$, and $X_z=\s _2^y
\s _1^x\in \Hom (x,z)$. Then
Thm.~\ref{th:Coxgr} and Prop.~\ref{pr:Stabrel} imply that
\begin{equation}
\begin{aligned}
  \Hom (x)\simeq \langle
  s_2,t,u_1,u_3\rangle /(&s_2^2,t^{m_{1,3}^x},u_1^2,u_3^2,
  (u_1s_2)^{m_{1,2}^x/3},\\
  &(s_2t^{-1}u_3t)^{m_{2,3}^x/3},(u_1u_3)^{m_{1,3}^z}),
\end{aligned}
  \label{eq:Homxcase2}
\end{equation}
where the inverse of the isomorphism is given by $s_2\mapsto \s _2^x$,
$t\mapsto \s _1^y\s _3^x$, $u_i\mapsto \s _1^y\s _2^z\s _i^z\s _2^y\s _1^x$
for $i\in \{1,3\}$.

If both $m_{1,2}^x$ and $m_{2,3}^x$ are even, then
\[\Hom (x)/(s_2)\simeq
\langle t,u_1,u_3\rangle
/(t^{m_{1,3}^x},u_1^2,u_3^2,(u_1u_3)^{m_{1,3}^z}).\]
Since $m_{1,3}^x\ge 2$, the group $\Hom (x)/(s_2)$ is not finite, which
is a contradiction to the finiteness of $\rsC $. By symmetry and by
Thm.~\ref{th:r2o3} we may assume that $m_{1,2}^x=3$. Then
Eq.~\eqref{eq:Homxcase2} tells that
\begin{equation}
\begin{aligned}
  \Hom (x)\simeq \langle
  s_2,u_3,t\rangle /(&s_2^2,u_3^2,t^{m_{1,3}^x},\\
  &(s_2t^{-1}u_3t)^{m_{2,3}^x/3},(s_2u_3)^{m_{1,3}^z}).
\end{aligned}
  \label{eq:Homxcase2-1}
\end{equation}
We apply Prop.~\ref{pr:Gmnp} to $\Hom (x)$. Since $m_{1,3}^x\ge 2$ and
$m_{1,3}^z\ge 2$, we conclude that
$\Hom (x)$ is finite if and only if
$(m_{1,3}^x,m_{2,3}^x)=(2,3)$ or $(m_{1,3}^x,m_{2,3}^x,m_{1,3}^z)=(3,3,2)$.
Hence $m_{2,3}^x=3$. Using
relation $m_{1,2}^x=3$, Lemma~\ref{le:Cartan-1}, (C2), and
Thm.~\ref{th:r2o3},
we obtain that
\begin{align*}
  \Cm ^x=
  \begin{pmatrix}
    2 & -1 & -a \\
    -1 & 2 & -1 \\
    -b & -1 & 2
  \end{pmatrix},
  \Cm ^y=
  \begin{pmatrix}
    2 & -1 & -a \\
    -1 & 2 & -1 \\
    -b & -1 & 2
  \end{pmatrix},
  \Cm ^z=
  \begin{pmatrix}
    2 & -1 & -c \\
    -1 & 2 & -1 \\
    -d & -1 & 2
  \end{pmatrix},
\end{align*}
where $a,b,c,d\in \ndN _0$. If $m_{1,3}^x=3$ and $m_{1,3}^z=2$ then $c=d=0$
and $a=b=1$ by Lemma~\ref{le:Cartan-1}. Then $a+1\not=c+1$, which is a
contradiction to Lemma~\ref{le:m=3}.
If $m_{1,3}^x=2$ then $a=b=0$,
and hence $c=d=0$ by Lemma~\ref{le:m=3} and (M2).
Then $\Cm ^x=\Cm ^y=\Cm ^z$ are Cartan matrices of type $A_3$. This proves the
theorem in Case~2.

\textit{Case 3:} Assume that $\rfl _3(x)=z$, $\rfl _3(y)=y$, and $\rfl _3(z)=x$.
For all
$a\in A$ and $i,j\in \{1,2,3\}$ with $i\not=j$ the relations
$(\rfl _i\rfl _j)^3(a)=a$ and $\rfl _i\rfl _j(a)\not=a$ hold,
and hence $3|m_{i,j}^a$ by (R4). In particular, $\cm ^a_{ij}<0$
for all $a\in A$ and $i,j\in \{1,2,3\}$ with $i\not=j$.

Calculate $\s :=\s _2^x\s _3^z\s _2^y\s _1^x\in \Hom (x)$. Since $\rsC $
is finite, $\s $ has finite order, and hence $\s =\id $ or $\tr \s <3$.
Direct calculation shows that $\s (\al _1)\in -\al _1+\ndZ \al _2+\ndZ \al _3$
and
\begin{align*}
  \tr \s =\cm ^x_{12}\cm ^x_{21}-\cm ^x_{12}\cm ^x_{23}\cm ^z_{31}
  &+(\cm ^x_{12}\cm ^y_{21}-1)(\cm ^x_{23}\cm ^z_{32}-1)\\
  &+\cm ^x_{13}\cm ^z_{31}
  -\cm ^x_{13}\cm ^y_{21}\cm ^z_{32}+\cm^y_{23}\cm ^z_{32}-2.
\end{align*}
By the above conclusion one obtains that $\tr \s \ge 3$ and $\s \not=\id $,
which is a contradiction. Thus there are no finite root systems of type $\cC $
in this case. This finishes the proof of the theorem.
\end{proof}

For the classification of root systems of type $\cC $ of rank higher than $3$
the following proposition will be useful.

\begin{propo}\label{pr:o=3,r=4}
  Let $\cC =\cC (I,A,(\rfl _i)_{i\in I},(\Cm ^a)_{a\in A})$
  be a connected Cartan scheme with
  $I=\{1,2,3,4\}$ and three objects $x,y,z$, and let
  $\rsC $ be a finite root system of type $\cC $.
  Assume that $\rfl _1(x)=y$, $\rfl _2(y)=z$, and
  $\cm ^a_{12}=\cm ^a_{21}=-1$ for all $a\in \{x,y,z\}$.
  Then $\rsC $ is standard.
\end{propo}

\begin{proof}
  \textit{Step 1. Both $\rsC |_{\{1,2,3\}}$ and $\rsC |_{\{1,2,4\}}$
  are standard.}

  Since $\rsC |_{\{1,2\}}$ is connected and irreducible, either
  $\rsC |_{\{1,2,3\}}=\rsC |_{\{1,2\}}\oplus \rsC |_{\{3\}}$,
  or $\rsC |_{\{1,2,3\}}$ is connected and irreducible.
  Then $\rsC |_{\{1,2,3\}}$ is standard --- in the
  first case, since $\rsC |_{\{1,2\}}$ is standard, and in the second case by
  Thm.~\ref{th:o=3,r=3}. Similarly, $\rsC |_{\{1,2,4\}}$ is standard.

  \textit{Step 2. If $\rfl _3(x)\not=x$ then $\rsC $ is standard.}

  By Thm.~\ref{th:r2o3}, $\rsC |_{\{2,3\}}$
  is connected and irreducible. Hence
  $\rsC |_{\{1,2,3\}}$ is connected and irreducible, and Thm.~\ref{th:o=3,r=3}
  implies that
  $\rfl _3(x)=y$ and $\cm ^a_{23}=\cm ^a_{32}=-1$ for all $a\in \{x,y,z\}$.
  Thus Step~1 gives that $\rsC |_{\{1,2,3\}}$, $\rsC
  |_{\{1,2,4\}}$, and $\rsC |_{\{3,2,4\}}$ are standard, hence $\rsC $ is
  standard, see Rem.~\ref{re:standard}.

  \textit{Step 3. $\rsC $ is standard.}

  If $\rfl _3\not=\id $ or $\rfl _4\not=\id $, then $\rsC $ is standard by
  Step~2. Assume now that $\rfl _3(a)=\rfl _4(a)=a$ for all $a\in A$.
  If $\rsC =\rsC |_{\{1,2\}}\oplus \rsC |_{\{3,4\}}$, then $\cm ^y_{ij}=0$
  for all
  $i\in \{1,2\}$, $j\in \{3,4\}$. Thus $\cm ^x_{jl}=\cm ^y_{jl}=\cm ^z_{jl}$
  for all $j\in \{3,4\}$ and $l\in I$ by Lemma~\ref{le:samerowfor0}.
  Since $\rsC |_{\{1,2,3\}}$ and $\rsC |_{\{1,2,4\}}$
  are standard by Step~1, it follows that 
  $\cm ^x_{jl}=\cm ^y_{jl}=\cm ^z_{jl}$ for all $j\in \{1,2\}$ and $l\in I$.
  Hence $\rsC $ is standard.

  Assume now that $\rsC |_{\{1,2,3\}}$ is irreducible. By Thm.~\ref{th:o=3,r=3}
  we may assume that $\cm ^a_{13}=\cm ^a_{31}=0$, $\cm ^a_{23}\cm ^a_{32}=2$,
  and $m_{2,3}^a=4$ for all $a\in \{x,y,z\}$.
  Then Thm.~\ref{th:Coxgr} and Prop.~\ref{pr:Stabrel} give that
  \begin{equation}
    \begin{aligned}
      \Hom (x)\simeq \langle &s_2,s_3,s_4,t_3,t_4,u_1,u_3,u_4\rangle /(
      s_2^2,s_3^2,s_4^2,t_3^2,t_4^2,u_1^2,\\
      &u_3^2,u_4^2,u_1s_2,s_3t_3,(s_4t_4)^{m_{1,4}^x/2},
      (s_2s_3)^4,\\
      &(s_2s_4)^{m_{2,4}^x},(s_3s_4)^{m_{3,4}^x},
      (t_3u_3)^2,(t_4u_4)^{m_{2,4}^y/2},\\
      &(t_3t_4)^{m_{3,4}^y},(u_1u_3)^2,(u_1u_4)^{m_{1,4}^z},
      (u_3u_4)^{m_{3,4}^z}
      ),
    \end{aligned}
    \label{eq:Homxr4}
  \end{equation}
  where the inverse of the isomorphism is given by the map
  $s_i\mapsto \s _i^x$ for all $i\in \{2,3,4\}$, $t_i\mapsto \s _1^y\s
  _i^y\s _1^x$ for all $i\in \{3,4\}$, and $u_i\mapsto \s _1^y\s
  _2^z\s _i^z\s _2^y\s _1^x$ for all $i\in \{1,3,4\}$.
  Therefore
  \begin{equation}
    \begin{aligned}
      \Hom &(x)\simeq \langle s_2,s_3,s_4,t_4,u_3,u_4\rangle /(
      s_2^2,s_3^2,s_4^2,t_4^2,u_3^2,u_4^2,\\
      &(s_4t_4)^{m_{1,4}^x/2}, (s_2s_3)^4,(s_2s_4)^{m_{2,4}^x},
      (s_3s_4)^{m_{3,4}^x},(s_3u_3)^2,\\
      &(t_4u_4)^{m_{2,4}^y/2},
      (s_3t_4)^{m_{3,4}^y},(s_2u_3)^2,
      (s_2u_4)^{m_{1,4}^z},(u_3u_4)^{m_{3,4}^z}
      ).
    \end{aligned}
    \label{eq:Homxr4-1}
  \end{equation}
  If $m_{1,4}^x>2$ and $m_{2,4}^y>2$, then $\Hom (x)$ is a Coxeter group
  without Coxeter relation between $s_2$ and $t_4$, which is a contradiction
  to the finiteness of $\rsC $.

  Assume first that $m_{1,4}^x=2$. Since $m_{1,3}^x=2$ as
  well, we obtain that $\cm ^x_{13}=\cm ^x_{14}=0$, and
  hence $m_{3,4}^x=m_{3,4}^y$. Then
  \begin{equation}
    \begin{aligned}
      \Hom &(x)\simeq \langle s_2,s_3,s_4,u_3,u_4\rangle /(
      s_2^2,s_3^2,s_4^2,u_3^2,u_4^2,(s_2s_3)^4,\\
      &(s_2s_4)^{m_{2,4}^x},
      (s_3s_4)^{m_{3,4}^x},(s_3u_3)^2,(s_4u_4)^{m_{2,4}^y/2},\\
      & (s_2u_3)^2,
      (s_2u_4)^{m_{1,4}^z},(u_3u_4)^{m_{3,4}^z}
      ).
    \end{aligned}
    \label{eq:Homxr4-2}
  \end{equation}
  If $m_{2,4}^y>2$, then there is no Coxeter relation between $s_4$ and $u_3$,
  which is a contradiction, and hence $m_{2,4}^y=2$.

  Similarly, if $m_{1,4}^x>2$ and $m_{2,4}^y=2$, then there is no Coxeter
  relation between $s_4$ and $u_3$.

  We are left with the case $m_{1,4}^x=m_{2,4}^y=2$.
  Then relations $\cm ^x_{14}=\cm ^y_{42}=0$ and
  Lemma~\ref{le:samerowfor0} imply that $\cm ^x_{42}=\cm ^y_{42}=0$, and hence
  $m_{2,4}^x=2$. Similarly, $\cm ^y_{24}=0$ implies that
  $\cm ^z_{41}=\cm ^y_{41}=\cm ^x_{41}=0$, hence $m_{1,4}^z=2$.
  Therefore
  \begin{equation}
    \begin{aligned}
      \Hom (x)\simeq \langle s_2,s_3,s_4,&u_3\rangle /(
      s_2^2,s_3^2,s_4^2,u_3^2,(s_2s_3)^4,(s_2s_4)^2,\\
      &(s_3s_4)^{m_{3,4}^x},(s_3u_3)^2,
      (s_2u_3)^2,(u_3s_4)^{m_{3,4}^z}
      ).
    \end{aligned}
    \label{eq:Homxr4-3}
  \end{equation}
  The finiteness of $\Hom (x)$ implies that $m_{3,4}^x,m_{3,4}^z\in \{2,3\}$.
  If $m_{3,4}^x=2$, then $m_{1,4}^x=m_{2,4}^x=0$ implies that $\rsC =\rsC
  |_{\{1,2,3\}}\oplus \rsC |_{\{4\}}$, and hence $\rsC $ is standard, since
  $\rsC |_{\{1,2,3\}}$ is standard. Further, $m_{3,4}^z=2$.
  Similarly, the assumption $m_{3,4}^z=2$ leads to the same result.
 
  Suppose now that $m_{3,4}^x=m_{3,4}^z=3$, that is,
  $\cm ^x_{34}=\cm ^x_{43}=\cm ^z_{34}=\cm ^z_{43}=-1$ by
  Lemma~\ref{le:Cartan-1}. Since $\cm^x_{13}=\cm ^x_{14}=0$,
  Lemma~\ref{le:samerowfor0} implies that $\cm ^y_{34}=\cm ^y_{43}=-1$, that
  is, $\Cm ^x=\Cm ^y=\Cm ^z$ is a Cartan matrix of type $F_4$. This proves the
  proposition.
\end{proof}

\begin{theor}\label{th:r>3o3}
  Let $\theta \in \ndN $ with $\theta \ge 4$, let
  $\cC =\cC (I,A,(\rfl _i)_{i\in I},(\Cm ^a)_{a\in A})$
  be a connected Cartan scheme with $I=\{1,2,\ldots ,\theta \}$
  and with $3$ objects, and let
  $\rsC $ be a finite irreducible root system of type $\cC $.
  Then $\rsC $ is standard, and the corresponding Cartan matrix $C$
  is indecomposable and of type $B_4$, $C_4$, $D_4$, or $F_4$.
  The object change diagram of $\rsC $ can be chosen to be
  \setlength{\unitlength}{1mm}
  \begin{align*}
    &
    \begin{picture}(26,6)(0,2)
      \put(3,3){\circle*{2}}
      \put(4.5,3.5){\line(1,0){7}}
      \put(4.5,2.5){\line(1,0){7}}
      \put(13,3){\circle*{2}}
      \put(14.5,3){\line(1,0){7}}
      \put(23,3){\circle*{2}}
      \put(8,5){\makebox[0pt]{$_{1,3}$}}
      \put(18,5){\makebox[0pt]{$_2$}}
    \end{picture}
    &&\text{if $C$ is of type $B_4$ or $C_4$,}\\&
    \begin{picture}(26,6)(0,2)
      \put(3,3){\circle*{2}}
      \put(4.5,3.7){\line(1,0){7}}
      \put(4.5,3){\line(1,0){7}}
      \put(4.5,2.3){\line(1,0){7}}
      \put(13,3){\circle*{2}}
      \put(14.5,3){\line(1,0){7}}
      \put(23,3){\circle*{2}}
      \put(8,5){\makebox[0pt]{$_{1,3,4}$}}
      \put(18,5){\makebox[0pt]{$_2$}}
    \end{picture}
    &&\text{if $C$ is of type $D_4$, and}\\&
    \begin{picture}(26,6)(0,2)
      \put(3,3){\circle*{2}}
      \put(4.5,3){\line(1,0){7}}
      \put(13,3){\circle*{2}}
      \put(14.5,3){\line(1,0){7}}
      \put(23,3){\circle*{2}}
      \put(8,5){\makebox[0pt]{$_1$}}
      \put(18,5){\makebox[0pt]{$_2$}}
    \end{picture}
    &&\text{if $C$ is of type $F_4$.}
  \end{align*}
\end{theor}

\begin{proof}
  As explained at the beginning of this section, we can choose $x,y,z\in A$
  and use a permutation of $I$ such that the object change diagram of
  $\rsC |_{\{1,2\}}$ is the one in Eq.~\eqref{eq:obchanger2}. Since $\rsC $ is
  irreducible, we can assume that $\rsC |_{\{1,2,3\}}$ is irreducible.
  Then Thm.~\ref{th:o=3,r=3} gives that
  $\rsC |_{\{1,2,3\}}$ is standard and $m^x_{1,2}\in \{2,3,4\}$.
  Further, since $\rfl _2\rfl _1(x)=z\not=y=\rfl _1\rfl _2(x)$ and
  $\rfl _1\rfl _2\rfl _1(x)=z=\rfl _2\rfl _1\rfl _2(x)$, (R4) yields
  that $m^x_{1,2}=3$. Hence $\cm ^a_{12}=\cm ^a_{21}=-1$ for all
  $a\in A$ by Lemma~\ref{le:Cartan-1}.
  By Prop.~\ref{pr:o=3,r=4} we obtain that $\rsC |_{\{1,2,i,j\}}$ is standard
  for all $i,j\in I\setminus \{1,2\}$ with $i\not=j$. Hence $\rsC $ is
  standard by Rem.~\ref{re:standard}.

  Finally, Thm.~\ref{th:o=3,r=3} and elementary argumentations
  with Dynkin diagrams imply that there is no connected
  irreducible standard root system, where the Cartan matrices are of type
  $A_4$ or of rank bigger than $4$,
  and that connected irreducible standard root systems of rank $4$ have
  the given object change diagrams. The rigorous proof is left to the reader.
\end{proof}

\section{Appendix}

In the following table we list all non-standard connected irreducible finite
Weyl groupoids with at most three objects. The last column contains the
stabilizer group of an object. It is a Coxeter group of the indicated type.

\vspace{11pt}
\begin{center}
  \begin{tabular}{| l c c c c c |}
    \hline
    $\Wg$ & $|A|$ & $|I|$ & $|\Wg|$ & $|R^a_+|$ & $\Hom(a)$\\
    \hline
    Thm.~\ref{th:o=2}(2) & 2 & 2 & 32 & 8 & $B_2$ \\
    Thm.~\ref{th:o=2}(2) & 2 & 2 & 48 & 12 & $G_2$ \\
    Thm.~\ref{th:o=2}(3) & 2 & 3 & 192 & 13 & $B_3$ \\
    Thm.~\ref{th:o=2}(4) & 2 & 3 & 192 & 13 & $B_3$ \\
    Thm.~\ref{th:r2o3}(3) & 3 & 2 & 36 & 6 & $A_1\times A_1$ \\
    Thm.~\ref{th:r2o3}(4) & 3 & 2 & 72 & 12 & $B_2$ \\
    Thm.~\ref{th:r2o3}(5) & 3 & 2 & 72 & 12 & $B_2$ \\
    Thm.~\ref{th:r2o3}(6) & 3 & 2 & 108 & 18 & $G_2$ \\
    Thm.~\ref{th:r2o3}(7) & 3 & 2 & 108 & 18 & $G_2$ \\
    \hline
  \end{tabular}
\end{center}


\providecommand{\bysame}{\leavevmode\hbox to3em{\hrulefill}\thinspace}
\providecommand{\MR}{\relax\ifhmode\unskip\space\fi MR }
\providecommand{\MRhref}[2]{%
  \href{http://www.ams.org/mathscinet-getitem?mr=#1}{#2}
}
\providecommand{\href}[2]{#2}

\end{document}